\documentclass[final,1p,times]{elsarticle}

\usepackage{amssymb}
\usepackage{amsmath}
\usepackage{algorithm}
\usepackage{algpseudocode}
\usepackage{bm}
\usepackage{booktabs}       
\usepackage{tabularx}
\usepackage{subcaption}

\usepackage{mathtools}
\usepackage{multirow}
\usepackage{bm}
\usepackage{makecell} 
\usepackage{setspace}

\usepackage{placeins}

\algnewcommand{\LineComment}[1]{\Statex {\(\triangleright\)} \textit{#1}}

\newcommand{\Ltwo}[1]{\|#1\|_{2}^2}

\usepackage{amsmath,amssymb,amsfonts}
\usepackage{amsthm}

\theoremstyle{plain} 
\newtheorem{theorem}{Theorem}[section]
\newtheorem{proposition}[theorem]{Proposition}

\newtheorem{corollary}[theorem]{Corollary}

\theoremstyle{definition} 
\newtheorem{definition}[theorem]{Definition}

\theoremstyle{remark} 
\newtheorem{remark}[theorem]{Remark}

\journal{Computers & Mathematics with Applications}

\begin{document}

\begin{frontmatter}




\title{Employing Deep Neural Operators for PDE control by decoupling training and optimization}

\affiliation[aff1]{organization={Department of Mathematics and Systems Analysis, Aalto University},
            city={Espoo},
            postcode={02150},
            country={Finland}}

\affiliation[aff2]{organization={DTU Management, Technical University of Denmark},
            city={Kgs. Lyngby},
            postcode={2800},
            country={Denmark}}

\author[aff1]{Oliver Lundqvist}
\ead{oliver.lundqvist@aalto.fi}

\author[aff1,aff2]{Fabricio Oliveira \corref{cor1}}
\ead{fabol@dtu.dk}

\cortext[cor1]{Corresponding author}

\begin{abstract}
Neural networks have been applied to control problems, typically by combining data, differential equation residuals, and objective costs in the training loss or by incorporating auxiliary architectural components. Instead, we propose a streamlined approach that decouples the control problem from the training process, rendering these additional layers of complexity unnecessary. In particular, our analysis and computational experiments demonstrate that a simple neural operator architecture, such as DeepONet, coupled with an unconstrained optimization routine, can solve tracking-type partial differential equation (PDE) constrained control problems with a single physics-informed training phase and a subsequent optimization phase. We achieve this by adding a penalty term to the cost function based on the differential equation residual to penalize deviations from the PDE constraint. This allows gradient computations with respect to the control using automatic differentiation through the trained neural operator within an iterative optimization routine, while satisfying the PDE constraints. Once trained, the same neural operator can be reused across different tracking targets without retraining. We benchmark our method on scalar elliptic (Poisson's equation), nonlinear transport (viscous Burgers' equation), and flow (Stokes equation) control problems. For the Poisson and Burgers problems, we compare against adjoint-based solvers: for the time-dependent Burgers problem, the approach achieves competitive accuracy with iteration times up to four times faster, while for the linear Poisson problem, the adjoint method retains superior accuracy, suggesting the approach is best suited to nonlinear and time-dependent settings. For the flow control problem, we verify the feasibility of the optimized control through a reference forward solver.
\end{abstract}


\begin{highlights}
\item \textbf{Decoupling}: A physics-informed neural operator is trained once and reused across multiple PDE-control objectives without retraining.
\item \textbf{Simplicity and surrogate use}: A plain DeepONet with a residual penalty solves tracking-type PDE-control problems without auxiliary networks or adjoint equations.
\item \textbf{Cheating directions}: Residual penalization provably suppresses non-physical optimization paths that exploit surrogate approximation errors.
\item \textbf{Computational advantage}: The surrogate approach achieves up to 4× faster solve times than adjoint methods on nonlinear time-dependent PDEs.
\end{highlights}


\end{frontmatter}



\section{INTRODUCTION}

A \textit{control problem} is an optimization problem in which the system dynamics are described by differential equations, either ordinary differential equations (ODEs) or partial differential equations (PDEs), that explicitly depend on a control input. A PDE-control problem is thereby an optimization problem where a PDE is the underlying constraint. These problems arise across a wide range of practical applications, such as optimization of flows \cite{ManservisiMenghini2016_CAMWA,CasasKunisch2019_SIAMJCO},
heat optimization \cite{ReinMohringDammKlar2020_DHN,KrugMehrmannSchmidt2020_DHN},
control of chemical reactions \cite{SeymenYucelKarasozen2014_DCR,CasasRyllTroeltzsch2018_RD},
and control of biological systems \cite{BongartiParkinsonWang2025_SIRRD,Yin2023optimal}.

Traditional numerical approaches to PDE-control problems are computationally intensive, making neural networks an attractive alternative. As such, recent studies have explored physics-informed neural networks (PINNs) and neural operators (NOs) for solving PDE-control and optimal control problems \cite{Molawi2023,verma2024neural,Alla2024,Tomasetto2024,Song2023admm}, particularly when the constraints are complex and nonlinear PDEs. Existing neural network models are typically trained for a single task, requiring retraining if problem parameters or cost functions change. For example, in a tracking problem, the tracked target may change, leading to a corresponding change in the cost function. Similarly, in an inverse setting where the goal is to estimate the source parameters that generated specific observations, any new measurement data would require retraining the model. It is worth noting that this reflects a critical practical limitation that prevents their deployment in real-world settings.

To address this limitation, we propose a method that fully decouples neural network training from the control problem, eliminating the need for retraining when cost functions or parameters change. Our approach trains a single physics-informed NO in advance, allowing it to solve control problems efficiently. We achieve that employing a direct optimization approach, integrating automatic differentiation with nonlinear optimization to efficiently compute gradients with respect to control decisions while enforcing differential equation constraints through a residual penalty. A schematic representation of our approach is shown in Figure \ref{fig:schematics}. To the best of our knowledge, and supported by our literature review, no prior work has shown that a plain physics-informed neural operator can be used as a direct surrogate for the PDE in an unconstrained optimization routine without auxiliary architectural components, without including the cost function in the training phase, or solving additional adjoint equations. 

In line with Occam’s razor, we advocate for simpler architectures, as our results demonstrate that such models are sufficient for solving a wide range of control problems without the need for elaborate or highly specialized designs. We show that solving a tracking PDE-control problem can be reduced to training a well-generalizable physics-informed NO, which precludes the need for specialized architectures and sophisticated numerical methods for solving PDEs and their adjoint equations. Consequently, our approach broadens the original scope of NOs beyond their typical use as solvers for differential equations, thereby showcasing a valuable discovery of their capabilities. We also show that the physics-informed NO can serve as an effective surrogate model for control problems, which is strongly supported by our numerical experiments. 

Our main contributions are as follows:

\begin{itemize}
    \item \textbf{Decoupled surrogate formulation.} We propose a two-stage framework in which a single physics-informed DeepONet is trained once to approximate the PDE solution operator and is then reused as a differentiable surrogate inside an unconstrained optimization routine to solve tracking-type PDE-control problems.
    \item \textbf{Residual-penalized objective for feasibility and stability.} We augment the surrogate tracking objective with a PDE-residual penalty and show how this term simultaneously enforces near-feasibility during the optimization trajectory and prevents the optimizer from exploiting surrogate approximation errors (i.e., suppresses non-physical ``cheating'' updates).
    \item \textbf{Reusability across objectives without retraining.} We demonstrate that the same trained neural operator can be applied to multiple tracking targets/cost functions for the same PDE class, eliminating retraining when the objective changes, thereby improving practical usability.
    \item \textbf{Empirical evaluation against adjoint baselines.} We benchmark the proposed method on scalar elliptic (Poisson's equation), nonlinear transport (viscous Burgers' equation), and a flow (Stokes equation) control tracking problems, and compare against classical adjoint-based solvers in terms of accuracy and computational runtime.
\end{itemize}
We limit our research to tracking-type PDE controls and to distributional controls, i.e., controls enforced directly in the PDE. We also consider the tracking targets to be reachable, i.e., that perfect tracking can be achieved to some discrepancy error. This problem type can also be interpreted as an inverse problem or control recovery. These assumptions also enable efficient numerical analysis of posterior errors and align with the PDE control literature (see, e.g.,  \cite{Yong2024DeepMRM,Hwang2022,Demo2023,Song2023admm}). We focus on the DeepONet architecture, though our approach is not limited to this architecture, leaving the exploration of alternatives for future work. Even though we focus on tracking costs and distributional control, our method shows promising results for other problem types, such as terminal cost or quadratic cost problems when a smoothing Tikhonov regularization is used in addition to the residual penalty, as we briefly show in \ref{appendix_d}. We note, however, that the mechanism behind this behaviour under additional regularization is not yet fully understood and thus is excluded from this study.

The rest of this paper is structured as follows: Section \ref{section_2} reviews related work of neural networks in optimal control problems in general. Section \ref{section_3} introduces our direct method as well as the theoretical justification of the method. Section \ref{section_4} describes the experimental setup on three PDE-control problems of tracking type. We optimize on two different cost functions for each PDE to demonstrate the reusability of the method. In section \ref{section_5} we present the results. Section \ref{section_6} concludes the study.

\section{RELATED WORK}\label{section_2}

In recent years, PINNs have been widely used to solve differential equations \cite{raissi2017physics,raissi2019physics,jagtap2022physics,kashefi2022physics,yang2021b}. An extension of PINNs for solving optimization problems was proposed by Lu et al. \cite{Lu2021}, where the authors used a PINN and a control network for topology optimization in inverse design problems. They trained a neural network on a common loss function consisting of the residual of the differential equation and the cost function of the topology optimization problem. Molawi et al. \cite{Molawi2023} extended the idea of a control neural network and a PINN to optimal control problems. Song et al. \cite{Song2023admm} split the control term into two variables, dividing the training loss into a differential equation term and a regularizing term describing the cost, and used an alternating direction method of multipliers for training the neural networks. Similarly, \cite{Demo2023,Alla2024,Tomasetto2024,barrystraume2025physics} have explored the use of two (or more) neural network models, PINNs and a control neural network, to solve optimal control problems. 

Another common approach is to employ PINNs together with the first-order necessary optimality conditions from Pontryagin’s Minimum Principle \cite{Pontryagin1956}. This yields a coupled system of differential equations for the state and adjoint variables, together with an optimality condition expressed in terms of the Hamiltonian. Yin et al. \cite{Yin2024AONN} used a direct adjoint looping approach to solve the optimal control problem by using two separate neural networks to solve it, one for the differential equation and one for the adjoint equation that arises from Pontryagin's Minimum Principle. Schiassi et al. \cite{schiassi2024pontryagin} also used Pontryagin's Minimum Principle and trained a PINN for solving the two-point boundary value problem for a quadratic cost function directly. Demo et al. \cite{Demo2023} used a PINN and a control network in sequence, where the output of the PINN was given as input to the control network, which solved for the adjoint equation and derived the optimal control. These methods showed promising results on their respective test problems, but are only capable of solving the particular control problem (i.e., with fixed cost function parameterization) for which the neural networks were trained.

Physics-informed NOs have also been applied to solve PDE-constrained control problems \cite{Yong2024Deep,Yong2024DeepMRM,Qi2024}. Hwang et al. \cite{Hwang2022} showed how a modified DeepONet can be used to solve control problems by first training a model and thereafter searching for the optimal control by using nonlinear programming and unconstrained optimization routines, where the gradients of the cost function were calculated by automatic differentiation with respect to the control. The approach in Hwang et al. \cite{Hwang2022} is similar to our proposed idea. However, the authors used an autoencoder to reconstruct the control function before re-feeding it to the NO model in the optimization routine's iterations. Further, for time-dependent PDEs,  in contrast to our approach, the authors used a transition network which predicts the next-time state from the current state and the control, similar to a time-integration scheme. Another recent approach is to include the Fréchet derivatives of the cost function in the training process \cite{cheng2025accelerating,luo2025efficient,yao2025difno}. This can improve gradient accuracy and accelerate optimization, but it also ties the trained model more closely to a particular objective functional and its derivatives, thereby reducing the degree of decoupling between surrogate training and the downstream control problem.

In contrast to the aforementioned approaches, we show that a physics-informed neural operator can solve PDE-constrained control problems by simply augmenting the cost function with a physics residual penalty, thereby avoiding additional architectural complexity and the more computationally expensive training processes associated with it.

\section{METHODOLOGY} \label{section_3}

\subsection{Deep Neural Operator}
Consider the following general partial differential equation:
\begin{equation}\label{diff}
\begin{aligned}
D(\mathbf{y}(\mathbf{x},t)) 
    &= f\!\bigl(\mathbf{y}(\mathbf{x},t), \mathbf{u}(\mathbf{x},t), \mathbf{x}, t\bigr),
    && (\mathbf{x},t)\in\Omega\times(0,T],\\
B(\mathbf{y}(\mathbf{x},t) )
    &= g(\mathbf{x},t),
    && (\mathbf{x},t)\in\partial\Omega\times(0,T],\\
\mathbf{y}(\mathbf{x},0) 
    &= h(\mathbf{x}),
    && \mathbf{x}\in\Omega.
\end{aligned}
\end{equation}
where $\mathbf{x}\in\Omega\subset\mathbb{R}^d$ denotes the spatial coordinate, and $t\in[0,T]$ is time. Here, $\mathbf{y}(\mathbf{x},t)\in\mathcal{Y}\subset\mathbb{R}^{n_y}$ is the \textit{state}, 
$\mathbf{u}(\mathbf{x},t) \in \mathcal{U} \subset \mathbb{R}^{n_u}$ is the \textit{control} function. The differential operator $D$ represents the system dynamics and
$f$ is a nonlinear function describing inputs and sources. The boundary operator $B$ encodes the boundary conditions on $\partial\Omega$. The function $g$ describes the boundary data, and $h$ prescribes the initial state. To ease the notation, we omit the arguments, for example, by denoting $\mathbf{u}(\mathbf{x},t):=\mathbf{u}$, only including them when the context requires. In addition, boldfaced symbols denote vector-valued quantities (e.g., $\mathbf{y}$, $\mathbf{u}$, $\mathbf{x}$), while non-boldfaced symbols denote scalar quantities, functions, or operators (e.g., $y$, $u$, $D$, $\mathcal{G}$).

We define $\mathbf{y}$ as the \textit{solution} of \eqref{diff}. We seek to find an operator $\mathcal{G}: \mathcal{U}\to \mathcal{Y}$ that solves the differential equation \eqref{diff} for any given function $\mathbf{u}$, where $\mathcal{U}$ and $\mathcal{Y}$ are appropriate Banach spaces. Hence, we write the \textit{solution operator} as
\begin{equation}\label{solution_operator}
\mathcal{G}(\mathbf{u})(\mathbf{x},t)=\mathbf{y}(\mathbf{x},t).
\end{equation}
A \textit{Neural Operator} approximates operator $\mathcal{G}$ with a neural network. Lu et al. \cite{lu2021learning} presented an architecture, the \textit{Deep Operator Network} (DeepONet), based on the universal approximation theorem for operators. Chen and Chen's \cite{chen1995universal} universal approximation theorem states that for any $\epsilon > 0$, there exist positive integers $n, p, m$, constants $c_i^k, \xi_{ij} \in \mathbb{R}$, points $\mathbf{x}_j \in \Omega$, and a continuous function $\sigma$, such that the operator $\mathcal{G}$ can be uniformly approximated by
\begin{equation}\label{deeponet_theorem}
    \left| \mathcal{G}(\mathbf{u})(\mathbf{x},t) - \sum_{k=1}^{p} \sum_{i=1}^{n} c_i^k \sigma\left( \sum_{j=1}^{m} \xi_{ij} \mathbf{u}(\mathbf{x}_j,t_j) \right) \sigma(\mathbf{w}_k \cdot (\mathbf{x}, t) + \zeta_k) \right| < \epsilon.
\end{equation}
The arguments of the continuous functions together with the functions themselves resemble two parallel neural networks, and therefore the theorem can be extended to learning solution operators with neural networks, especially in the form of the DeepONet \cite{lu2021learning}. The basic architecture of the DeepONet consists of two parallel networks called \textit{branch} and \textit{trunk}. The trunk network receives discrete space–time coordinates $(\mathbf{x}_i,t_i)$, where the single index $i$ enumerates the discretization points in both space and time. Respectively,  the branch network receives samples of the input function $\mathbf{u}$ at $m$ fixed sensor locations $S=\{\mathbf{s}_j\}_{j=1}^{m}$ where $\mathbf{s}_j=(\mathbf{x}_j,t_j)\in\Omega$. The stacked sensor values at the input are denoted as $\mathbf{u}_S=[\mathbf{u}(\mathbf{s}_1),\dots,\mathbf{u}(\mathbf{s}_m)]\in\mathbb{R}^{n_u\times m}$. At the end of the DeepONet, the outputs of the trunk and branch networks are combined with a dot product, which can be compactly written as
\begin{equation}\label{deeponet_approx}
    \mathcal{G}_{\bm{\theta}}(\mathbf{u}_S)(\mathbf{x}_i,t_i) = \sum_{k=1}^{p} \underbrace{b_k\big(\mathbf{u}(\mathbf{s}_1), \dots, \mathbf{u}(\mathbf{s}_{m})\big)}_{\text{Branch Network}} \cdot \underbrace{\tau_k(\mathbf{x}_i,t_i)}_{\text{Trunk Network}} +\; \mathbf{b}_0,
\end{equation}
where $b_k$ and $\tau_k$ are the neural outputs of the branch and trunk networks, respectively, $p$ is the number of neurons in the last layer of the branch and trunk network, and $\mathbf{b}_0 \in \mathbb{R}^{n_y}$ is an optional bias vector. Both the trunk and branch networks are fully connected networks in their simplest form. A fully connected network (FCN) with $L$ layers is a composition of affine maps and elementwise nonlinearities. Let $\sigma$ be the activation function, then we define the FCN as 
\begin{equation}\label{eq:MLP_general}
\begin{aligned}
\mathrm{FCN}_{\theta}(\mathbf{z})
&:= \mathbf{h}^{(L)}(\mathbf{z}), \qquad \mathbf{h}^{(0)}(\mathbf{z}) := \mathbf{z},\\
\mathbf{h}^{(\ell)}(\mathbf{z})
&:= \sigma\!\left(\mathbf{W}^{(\ell)}\mathbf{h}^{(\ell-1)}(\mathbf{z})+\mathbf{c}^{(\ell)}\right),
\qquad \ell=1,\dots,L-1,\\
\mathbf{h}^{(L)}(\mathbf{z})
&:= \mathbf{W}^{(L)}\mathbf{h}^{(L-1)}(\mathbf{z})+\mathbf{c}^{(L)}.
\end{aligned}
\end{equation}%
Here, $\theta := \{\mathbf{W}^{(\ell)},\mathbf{c}^{(\ell)}\}_{\ell=1}^{L}$ denotes all trainable parameters. The last layer is typically kept linear. Therefore we define the NO $\mathcal{G}_{\theta}:\mathbb{R}^{n_u\times m}\times \Omega\times [0,T]\rightarrow\mathbb{R}^{n_y}$ with parameters $\theta$ and at a discrete space–time point $(\mathbf{x}_i,t_i)$ as
\begin{equation}\label{NO_solution_operator}
\mathcal{G}_\theta(\mathbf{u}_S)(\mathbf{x}_i,t_i)=\mathbf{y}_{\theta,i},
\end{equation}
where $\mathbf{y}_\theta\approx \mathbf{y}$. In practice, multiple coordinate points $\{(\mathbf{x}_i,t_i)\}$ can be input as a batch for efficient evaluation, but we omit this in the notation for clarity. For more details, we refer the reader to \cite{lu2021learning,wang2021learning,wang2021understanding}. 

Training the NO can be performed in a purely data-driven fashion using pairs of input functions and corresponding solutions (e.g., minimizing \ $\|\mathbf{y}_\theta-\mathbf{y}\|_2^2$), or in a \textit{physics-informed} way. In the latter, the loss function incorporates the residual of the governing differential equation, as well as boundary and initial conditions. This enforces physical consistency by penalizing violations of $D(\mathbf{y}_\theta)-f=0$, $B(\mathbf{y}_\theta)-g=0$ and $\mathbf{y}_\theta(\mathbf{x},0)-h=0$, enabling training even without data, while embedding the dynamics directly into the learning process. We define the discrete \textit{differential residual} function  as
\begin{equation}\label{residual}
\mathcal{R}_{\theta,i} :=\mathcal{R}_{\theta,i}(\mathbf{y}_{\theta,i},\mathbf{u}_{i},\mathbf{x}_i,t_i)= D(\mathbf{y}_{\theta,i}) - f(\mathbf{y}_{\theta,i}, \mathbf{u}_{i}, \mathbf{x}_i,t_i),
\end{equation}
and the \textit{constraint residuals}, which includes the boundary 
\begin{equation}\label{residual_boundary}
\mathcal{R}_{\theta,i}^B :=\mathcal{R}_{\theta,i}^B(\mathbf{y}_{\theta,i},\mathbf{x}_i,t_i)= B(\mathbf{y}_{\theta,i}) -g(\mathbf{x}_i,t_i), 
\end{equation}
and initial conditions
\begin{equation}\label{residual_initial}
\mathcal{R}_{\theta,i}^{t_0} :=\mathcal{R}_{\theta,i}^{t_0}(\mathbf{y}_{\theta,i},\mathbf{x}_i) =\mathbf{y}_{\theta,i} -h(\mathbf{x}_i). 
\end{equation}
A central element in this way of training the NO is the differential operation $D(\mathbf{y}_{\theta})$, which can be easily computed with automatic differentiation at the discrete points. We denote with the subscript $\theta$ that the residuals in equations \eqref{residual}-\eqref{residual_initial} are evaluated by automatic differentiation through the neural operator. We summarize the process of training a physics-informed NO by solving an optimization problem
\begin{equation}\label{physics_informed_training}
\min_{\theta} \; L(\theta):=\frac{1}{|I|}\sum_{i\in I} \|\mathcal{R}_{\theta,i}\|^2+ \frac{1}{|\partial I|}\sum_{j\in \partial I}\|\mathcal{R}^B_{\theta,j}\|^2 + \frac{1}{|I_0|}\sum_{k\in I_0}\|\mathcal{R}^{t_0}_{\theta,k}\|^2,
\end{equation}
where $I$ denotes the set of inner collocation points, $\partial I$ boundary collocation points and $I_0$ the collocation points when $t=0$. For more technical details, implementations, and variants of physics-informed training, we refer to \cite{raissi2017physics,raissi2019physics,wang2021learning}.

\subsection{Physics Informed Neural Operators for Direct Methods}
A typical continuous PDE-constrained tracking problem can be defined as
\begin{equation}\label{eq:pde_tracking_cont}
\begin{aligned}
\min_{\mathbf{u}\in\mathcal{U}}\quad 
J(\mathbf{y},\mathbf{u})
:=\;& \int_{0}^{T}\int_{\Omega}
\|\mathbf{y}(\mathbf{x},t)-\mathbf{y}_d(\mathbf{x},t)\|_2^2\,d\mathbf{x}\,dt
\;+\;\lambda\|\mathbf{u}\|_{\mathcal{U}}^2 \\
\text{s.t.}\quad 
& D(\mathbf{y}) = f(\mathbf{y},\mathbf{u},\mathbf{x},t) \\
& B(\mathbf{y}) = g(\mathbf{x},t) \\
& \mathbf{y}(\mathbf{x},0) = h(\mathbf{x}),
\end{aligned}
\end{equation}
where the objective is to find a control function $\mathbf{u}$ in some appropriate Hilbert space $\mathcal{U}$ (e.g., $L_2$) that minimizes the objective $J$. Here, the quadratic term $\|u\|_\mathcal{U}^2$ provides Tikhonov regularization, which stabilizes the otherwise ill-posed tracking/inverse problem and yields a well-posed optimization problem. The constraint equation $D(\mathbf{y}) 
= f(\mathbf{y}, \mathbf{u}, \mathbf{x}, t)$ is a set of differential equations with the boundary and initial conditions $B(\mathbf{y})=g(\mathbf{x},t)$ and $\mathbf{y}(\mathbf{x},0)=h(\mathbf{x})$.

For solving (\ref{eq:pde_tracking_cont}), one can resort to the \textit{indirect methods}, also known as \textit{optimize-then-discretize}, and use the \textit{adjoint method} to derive the differential equation of adjoint variables, which are used for computing the gradients of the cost function \cite{Hinze2009OptimizationPDEConstraints,lions1971optimal}. For this, one formulates a Lagrangian functional $\mathcal{L}$ that adjoins the PDE constraints to the objective function using an \textit{adjoint variable} (or co-state), denoted as $\mathbf{p}(\mathbf{x},t)$. Let $\mathbf{e}(\mathbf{y},\mathbf{u}):=D(\mathbf{y})- f (\mathbf{y},\mathbf{u},\mathbf{x},t)=0$ be the exact PDE in implicit form. Then the Lagrangian functional without terminal cost can be stated as
\begin{equation}\label{eq:lagrangian}
\mathcal{L}(\mathbf{y}, \mathbf{u}, \mathbf{p}) = J(\mathbf{y}, \mathbf{u}) + \int_{0}^{T}\int_{\Omega} \mathbf{p}^{\top} \mathbf{e}(\mathbf{y},\mathbf{u})\,d\mathbf{x}\,dt.
\end{equation}
By using the Fréchet derivative of the Lagrangian with respect to $\mathbf{y}$ and setting it to zero, one yields the \textit{adjoint equation} for Dirichlet boundary conditions:
\begin{equation}\label{eq:adjoint}
\begin{aligned}
\nabla_\mathbf{y}J(\mathbf{y},\mathbf{u}) + (\nabla_{\mathbf{y}} \mathbf{e})^* \mathbf{p}&= \mathbf{0}, \quad &&(\mathbf{x},t)\in\Omega\times[0,T), \\
B^*(\mathbf{p}) &= \mathbf{0}, \quad &&(\mathbf{x},t)\in\partial\Omega\times[0,T), \\
\mathbf{p}(\mathbf{x}, T) &= \mathbf{0}, \quad &&\mathbf{x}\in\Omega,
\end{aligned}
\end{equation}
where $(\nabla_{\mathbf{y}} \mathbf{e})^*$ is the adjoint of the Jacobian of $\mathbf{e}$ with respect to $\mathbf{y}$. For a tracking problem the term $\nabla_\mathbf{y}J(\mathbf{y},\mathbf{u})$ reduces down to $2(\mathbf{y}-\mathbf{y}_d)$. By solving the state equation forward in time and the adjoint equation backwards in time, the gradient of the objective functional can be evaluated purely in terms of $\mathbf{y}$ and $\mathbf{p}$. Hence, for a given $\mathbf{u}$,  
\begin{equation}\label{eq:adjoint_gradient}
\nabla_{\mathbf{u}} \mathcal{L} = \nabla_{\mathbf{u}} J(\mathbf{y},\mathbf{u)} + (\nabla_{\mathbf{u}} \mathbf{e})^* \mathbf{p}.
\end{equation}
The gradient \eqref{eq:adjoint_gradient} can then be used in any gradient descent algorithm for solving the PDE-control problem, by minimizing the tracking cost. For a comprehensive derivation of the Lagrangian formalism and further theoretical details on adjoint-based PDE-constrained optimization, we refer the reader to \cite{Hinze2009OptimizationPDEConstraints, Troeltzsch2010OptimalControlPDE}.

Another approach is the \textit{direct method} or \textit {discretize-then-optimize}, which discretizes the optimal control problem and treats it as a nonlinear programming (NLP) problem \cite{Hinze2009OptimizationPDEConstraints}. This can be achieved by first discretizing the space–time domain into grid points indexed by $i=0,\dots, I$. The discrete approximations of the functions are denoted as $\mathbf{y}_i \approx \mathbf{y}(\mathbf{x}_i,t_i)$ and $\mathbf{u}_i \approx \mathbf{u}(\mathbf{x}_i,t_i)$ at the $i$-th grid point. Further, we note the stacked discretized control and state vectors as $\mathbf{u}_h=[\mathbf{u}_1,\dots \mathbf{u}_I]$ and $\mathbf{y}_h=[\mathbf{y}_1,\dots,\mathbf{y}_I]$. This yields the discrete optimization problem $\min \bar{J}(\mathbf{y}_h,\mathbf{u}_h)$ as
\begin{equation}\label{eq:pde_tracking_discrete}
\begin{aligned}
\min_{\mathbf{u}_h}\quad 
& \bar{J}(\mathbf{y}_h,\mathbf{u}_h)
:= \sum_{i=1}^{I}\|\mathbf{y}_i-\mathbf{y}_{d,i}\|_2^2\,w_i
\;+\; \lambda\,\|\mathbf{u}_h\|_{\mathcal{U}_h}^2 .\\
\text{s.t.}\quad 
& D(\mathbf{y}_i) = f(\mathbf{y}_i,\mathbf{u}_i,\mathbf{x}_i,t_i), \quad \forall(\mathbf{x}_i,t_i)\in \Omega\times[0,T],\\
& B(\mathbf{y}_j) = g(\mathbf{x}_j,t_j), \; \quad \qquad  \forall(\mathbf{x}_j,t_j)\in \partial\Omega.\\
&\mathbf{y}(\mathbf{x}_k,0)=h(\mathbf{x}), \qquad \qquad \forall \mathbf{x}_k\in \Omega,
\end{aligned}
\end{equation}
where $w_i$ are quadrature weights from the chosen integration rule, such as Gaussian quadrature or trapezoidal integration. We seek to find the control $\mathbf{u}_h\in\mathcal{U}_h$ where $\mathcal{U}_h\subset\mathcal{U}$ is a finite-dimensional subspace of the Hilbert space $\mathcal{U}$. The single index $i$ implicitly enumerates both space and time discretization points. Equation (\ref{eq:pde_tracking_discrete}) is now an NLP problem with only equality constraints. 

In general, solving NLP problems requires the computation of gradients of the cost function $\bar{J}(\mathbf{u}_h)$ with respect to $\mathbf{u}_h$, i.e., $\nabla_\mathbf{u}\bar{J}(\mathbf{u}_h)$. Thus, solving a discretized optimal control problem directly is not trivial as it requires the computation of the Jacobian of the discretized state vector, i.e., the \textit{sensitivities} $\nabla_\mathbf{u}\mathbf{y}_h$. 

An attractive choice to solve \eqref{eq:pde_tracking_discrete} is to use a NO as a surrogate model for efficiently solving the PDEs. That is, we replace the differential equations in the constraints of \eqref{eq:pde_tracking_discrete} with our pre-trained NO. To ensure that the integrals in the cost function are evaluated without the need for interpolation, the coordinate discretization points $(\mathbf{x}_i,t_i)$ are chosen to coincide with the sensor locations $S$, i.e., $\mathbf{u}_h=\mathbf{u}_S$ and $\mathbf{y}_h=\mathbf{y}_S$. This leads to the following formulation:
\begin{equation}\label{eq:pde_tracking_no_reduced}
\min_{\mathbf{u}_S}\;
\bar{J}_{NO}(\mathbf{u}_S)
:=
\sum_{j=1}^{m}
\left\|\mathcal{G}_\theta(\mathbf{u}_S)(\mathbf{s}_j)-\mathbf{y}_{d}(\mathbf{s}_j)\right\|_2^2\,w_j
\;+\;
\lambda\,\|\mathbf{u}_S\|_{\mathcal{U}_h}^2,
\end{equation}
where $\mathbf{s}_j=(\mathbf{x}_j,t_j)\in\Omega$ denotes the sensor locations.
Since the NO $\mathcal{G_\theta}$ is at its core a neural network, its output can be differentiated with respect to the input; thus, an optimization routine for NLP problems, such as gradient descent, can be employed to solve the problem. For that, we require two things. First, the NO must be trained, generalizable, and expressive enough to represent a sufficiently large near-feasible set, i.e., a region where the PDE residual remains small along the optimization trajectory. This can be ensured in practice by making the number of functions in the NO training set sufficiently large, such that it contains a variety of functions. Second, since we want to optimize on the control (input) $\mathbf{u}$ and repeatedly update $\mathbf{u}$, the NO can produce solutions $\mathbf{y}$ that do not approximate the solution to the differential equation. To address this, the cost function must be properly penalized with the residual of the differential equation, ensuring that the control remains in the feasible space satisfying the differential equation. Without proper penalization driven by this residual, the computed gradients are noisy and can exploit surrogate errors, and thereby fail to find a feasible solution as we discuss in section \ref{section3.4} and demonstrate in section \ref{section_5}. Thus, we penalize the cost function \eqref{eq:pde_tracking_no_reduced} with the quadratic differential residual \eqref{residual}, obtaining%
\begin{equation}\label{eq:pde_tracking_no_penalized_reduced}
\begin{aligned}
\min_{\mathbf{u}_S}\;
\bar{J}_{\mu}(\mathbf{u}_S)
:=&
\sum_{j=1}^{m}
\left\|\mathcal{G}_\theta(\mathbf{u}_S)(s_j)-\mathbf{y}_{d}(s_j)\right\|_2^2\,w_j
\;\\
+&\lambda\,\|\mathbf{u}_S\|_{\mathcal{U}_h}^2 +
\mu\|\mathcal{R}_\theta(\mathbf{u}_S)\|^2_2, 
\end{aligned}
\end{equation}
where $\mu>0$ is a \textit{penalty} factor, and we have defined the residuals evaluated through the neural operator as
\begin{equation}\label{eq:residual_NO}
\mathcal{R}_\theta(\mathbf{u}_S)
:= [\mathcal{R}_\theta(\mathbf{u}(\mathbf{s}_1)),\dots,\mathcal{R}_\theta(\mathbf{u}(\mathbf{s}_m))]
\end{equation}
by using automatic differentiation and calculating the respective residuals with \eqref{residual}.
Notice that problem \eqref{eq:pde_tracking_no_penalized_reduced} only considers the physics residual $\mathcal{R}_\theta$. Similarly, the boundary residual $\mathcal{R}_\theta^B$ and initial condition residual $\mathcal{R}_\theta^{t_0}$ could be added to the cost function as a penalty, but our experiments showed it provided no further benefit. Therefore, we omit them as they are not necessary in practice.

We note that employing regularization, such as the Tikhonov regularization $\Ltwo{\mathbf{u}}$ or the $H^1$-seminorm regularization $\Ltwo{\nabla \mathbf{u}}$ is optional, but can provide smoother solutions and reduce noise. Since the NO $\mathcal{G}_\theta$ is a composition of differentiable operations (e.g., linear layers and activation functions), the penalized objective functional $\bar{J}_\mu(\mathbf{u}_S)$ is differentiable with respect to the control input $\mathbf{u}_S$. Given a sufficiently large penalty factor $\mu$, it enforces the optimization method to remain within the domain where the residuals of the differential equation remain small. Moreover, since the NO was trained physics-informed, it is capable of predicting results that do not violate the dynamical constraints.

\begin{figure}[t]  
    \centering
    \includegraphics[width=0.9\linewidth]{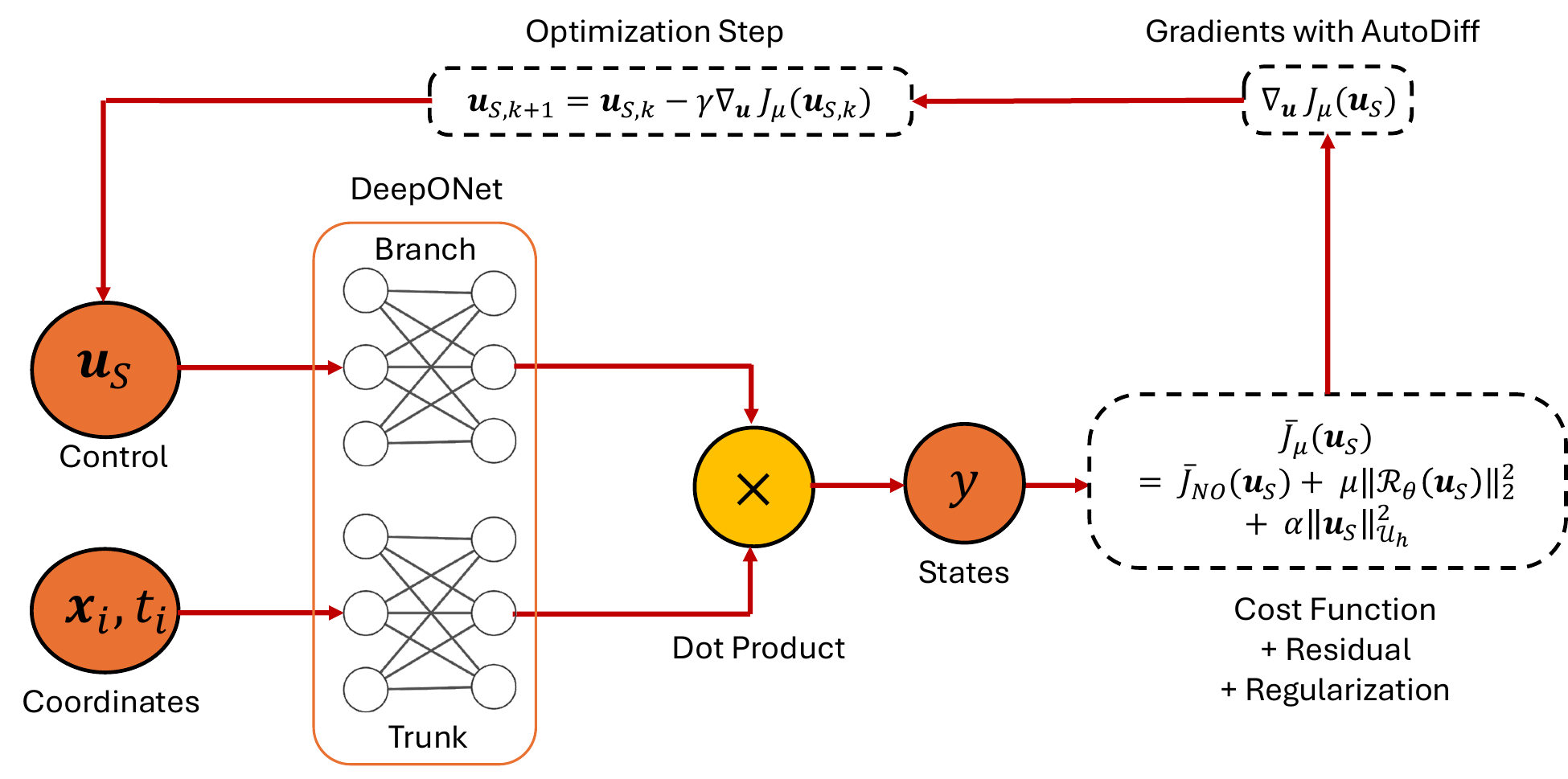}  
    \caption{Schematic of the proposed framework. 
    A NO (DeepONet) maps the control input $u$ and 
    space--time coordinates $(\mathbf{x},t)$ into the system state $y$. 
    The state is used to evaluate a penalized objective $\bar{J}_\mu(u)$, 
    which combines the original cost function with physics residuals 
    and regularization terms. Gradients $\nabla_\mathbf{u}\bar{J}_\mu(\mathbf{u}_S)$ 
    are computed via automatic differentiation, and the control is updated 
    through an optimization step. This process iterates until convergence, 
    yielding the optimal control.}
    \label{fig:schematics}
\end{figure}

Formulated as a standard NLP, we seek to find a control $\mathbf{u}_S$ by fixing the network parameters $\theta$ and minimizing $\bar{J}_\mu$ via an iterative gradient-based optimization method. In each iteration $k$, the gradient of the objective with respect to the control, $\nabla_u \bar{J}_\mu(\mathbf{u}_{S,k})$, is computed using automatic differentiation and backpropagated through the NO. The control is then updated using a descent step, such as:
\begin{equation}\label{eq:gradient_descent}
\mathbf{u}_{S,{k+1}} = \mathbf{u}_{S,k} - \gamma \nabla_u \bar{J}_\mu(\mathbf{u}_{S,k}),
\end{equation}
where $\gamma > 0$ is the step size. This approach effectively treats the pre-trained physics-informed neural operator as a differentiable surrogate for the differential equations, allowing the use of efficient off-the-shelf unconstrained optimizers (e.g., L-BFGS or Adam) to solve the control problem without repeatedly querying a numerical PDE solver or solving the adjoint equations for obtaining gradients. The method is schematically presented in Figure \ref{fig:schematics}. With a slight abuse of terminology, we refer to our method hereinafter as the \textit{penalty method} (see Remark \ref{remark:penalty}).

\subsection{Stationary conditions}

Next, we provide analytical justification for employing residual penalization and state the corresponding stationarity conditions. For notational convenience, we absorb boundary and initial conditions into the residual operator and work entirely in a reduced formulation. We consider a pre-trained NO $\mathcal{G}_\theta$ and define the NO-based constrained problem as
\begin{equation}\label{eq:pde_simple}
\begin{aligned}
\min_{\mathbf{u}_S}\quad 
& \bar{J}_{NO}(\mathbf{u}_S)\\
\text{s.t.}\quad 
& \mathcal{R}_\theta(\mathbf{u}_S)=\mathbf{0},
\end{aligned}
\end{equation}
With this notation, the relaxed NO-based unconstrained problem is:
\begin{equation}\label{eq:pde_penalized_reduced}
\min_{\mathbf{u}_S}\;
\bar{J}_{\mu}(\mathbf{u}_S)
:=\bar{J}_{\mathrm{NO}}(\mathbf{u}_S) +\mu\|\mathcal{R}_\theta(\mathbf{u}_S)\|^2_2,
\end{equation}
where $\mu>0$ is the fixed penalty parameter. Problem \eqref{eq:pde_penalized_reduced} is now an unconstrained optimization problem. Since $\mathcal{G}_\theta$, $\bar{J}_{\mathrm{NO}}$ and $\mathcal{R}_\theta$ are differentiable with respect to $\mathbf{u}_S$, it allows us to pose stationarity conditions for the penalized problem \eqref{eq:pde_penalized_reduced} that are analogous to the Karush-Kuhn-Tucker (KKT) optimality conditions, yielding the following proposition.

\begin{proposition}[stationarity condition]\label{prop:penalty_kkt_like}
Let $\mathbf{u}_\mu$ be a stationary point of the penalized objective
\eqref{eq:pde_penalized_reduced}. Then $\mathbf{u}_\mu$ satisfies the identity
\begin{equation}\label{eq:grad_decomposition}
\nabla_{\mathbf{u}} \bar{J}_{\mu}(\mathbf{u}_\mu)
=
\nabla_{\mathbf{u}} \bar{J}_{NO}(\mathbf{u}_\mu)
+2\mu\,\Big({\nabla_{\mathbf{u}}\mathcal{R}_\theta(\mathbf{u}_\mu})\Big)^\top \mathcal{R}_\theta(\mathbf{u}_\mu)
=\mathbf{0},
\end{equation}
with the implicit multiplier
\begin{equation}\label{eq:implicit_multiplier}
\boldsymbol{\nu}_\mu := 2\mu\,\mathcal{R}_\theta(\mathbf{u}_\mu),
\end{equation}
which yields the stationarity condition
\begin{equation}\label{eq:kkt_like_stationarity}
\nabla_{\mathbf{u}} \bar{J}_{NO}(\mathbf{u}_\mu) + ({\nabla_{\mathbf{u}}\mathcal{R}_\theta(\mathbf{u}_\mu}))^\top \boldsymbol{\nu}_\mu = \mathbf{0}.
\end{equation}
\end{proposition}
\begin{proof}
The proof follows directly from the linearity of the gradient, applying the chain rule on the penalty, and regrouping and renaming the terms.
\end{proof}
The implication of Proposition \ref{prop:penalty_kkt_like} is that at convergence, the penalized method produces a point that is stationary for the penalized NLP \eqref{eq:pde_penalized_reduced} and approximately feasible when the PDE-residual $\mathcal{R}_\theta$ is small. Further, the tracking gradient $\nabla_\mathbf{u}\bar{J}_{NO}$ is balanced by a physics-consistency correction term $2\mu(\nabla_{\mathbf{u}}\mathcal{R}_\theta)^\top \mathcal{R}_\theta$. In practice, this means that the correction term $2\mu(\nabla_{\mathbf{u}}\mathcal{R}_\theta)^\top \mathcal{R}_\theta$ prevents the optimizer from exploiting the NO errors and high frequency directions that would reduce the tracking loss while violating the PDE-constraints. 

\begin{remark}\label{remark:penalty}
In contrast to classical penalty methods, where the sequence of solutions $\{\mathbf{u}_\mu\}$ to the penalized problem \eqref{eq:pde_penalized_reduced} converges to the solution of constrained problem \eqref{eq:pde_simple} as $\mu\to\infty$ (see, for example, \cite{bazaraa2006nonlinear}), we do not increase the penalty factor and keep it constant during the iterations. In other words, we solve a regularized problem instead and obtain stationary points that satisfy the conditions of Proposition \ref{prop:penalty_kkt_like}. Thus, 
 albeit being a slight abuse of terminology regarding the classical penalty method, our ``penalty'' method does share similarities by turning the constrained problem into an unconstrained one and using the residual as guiding the solution or to remain in the feasible domain, similarly to the classical penalty method.
\end{remark}

The gradient of the penalized reduced objective admits a factorization that implicitly encodes the adjoint equation and gradient of the Lagrangian, which motivates the following proposition.

\begin{proposition}[Implicit formulation of the adjoint equation and gradient]\label{prop:implicit_adjoint}
Let the NO-evaluated residual by composition be $\mathcal{R}_\theta(\mathbf{u}_S):=\mathbf{e}\big(\mathbf{y}_\theta,\,\mathbf{u}_S\big)$, where $\mathbf{y}_\theta:= \mathcal{G}_\theta(\mathbf{u}_S)$. Then, for every $\mathbf{u}_S$ (where the derivatives exist), the gradient of the penalized reduced objective, admits the exact factorization 

\begin{equation}\label{eq:implicit_adjoint_grad}
\begin{aligned}
\nabla_\mathbf{u}\bar{J}_\mu=\nabla_{\mathbf{u}}\mathcal{G}_\theta(\mathbf{u}_S)^\top \Big[ \nabla_{\mathbf{y}} J(\mathbf{y}_{\theta}, \mathbf{u}_S) + \nabla_{\mathbf{y}} \mathbf{e}(\mathbf{y}_{\theta}, \mathbf{u}_S)^\top \boldsymbol{\nu}_\mu \Big] \\
+ \Big[ \nabla_{\mathbf{u}} J(\mathbf{y}_{\theta}, \mathbf{u}_S) + \nabla_{\mathbf{u}} \mathbf{e}(\mathbf{y}_{\theta}, \mathbf{u}_S)^\top \boldsymbol{\nu}_\mu \Big] ,
\end{aligned}
\end{equation}
where the terms in the brackets correspond exactly to the classical adjoint equation and gradient of the Lagrangian, respectively, evaluated at the surrogate state with the implicit multiplier $\boldsymbol{\nu}_\mu := 2\mu\mathcal{R}_\theta(\mathbf{u}_S)$.
\end{proposition}

\begin{proof}
Let $\mathbf{u}_S$ be Fréchet differentiable with the state $\mathbf{y}_\theta = \mathcal{G}_\theta(\mathbf{u}_S)$. Then the gradient terms of Proposition \ref{prop:penalty_kkt_like} can be written as
\begin{equation}
\nabla_\mathbf{u}\bar{J}_{NO}=(\nabla_\mathbf{y}J(\mathbf{y}_\theta,\mathbf{u}_S))^\top\nabla_{\mathbf{u}}\mathcal{G}_\theta(\mathbf{u}_S)+ \nabla_{\mathbf{u}}J(\mathbf{y}_\theta,\mathbf{u}_S),
\end{equation}
and
\begin{equation}
({\nabla_{\mathbf{u}}\mathcal{R}_\theta(\mathbf{u}_S}))^\top \boldsymbol{\nu}_\mu = (\nabla_{\mathbf{u}} \mathbf{e}(\mathcal{G}_\theta(\mathbf{u}_S),\mathbf{u}_S))^\top \boldsymbol{\nu}_\mu=
(\nabla_{\mathbf{y}} \mathbf{e}(\mathbf{y}_\theta, \mathbf{u}_S)\nabla_\mathbf{u}\mathcal{G}_\theta(\mathbf{u}_S)+ \nabla_\mathbf{u}\mathbf{e}(\mathbf{y}_\theta, \mathbf{u}_S))^\top\boldsymbol{\nu}_\mu.
\end{equation} 
Then factorizing with respect to $\nabla_\mathbf{u}\mathcal{G}_\theta(\mathbf{u}_S)$ and regrouping the terms, gives \eqref{eq:implicit_adjoint_grad}.
\end{proof}

The usefulness of Proposition~\ref{prop:implicit_adjoint} relates to the structure it exposes. Equation~\eqref{eq:implicit_adjoint_grad} shows that the reduced gradient is obtained by evaluating the adjoint equation residual \eqref{eq:adjoint} and Lagrangian gradient \eqref{eq:adjoint_gradient} at the surrogate state $\mathbf{y}_\theta=\mathcal{G}_\theta(\mathbf{u}_S)$ and mapping the resulting state-residual back to the control variables through the sensitivity $\nabla_{\mathbf{u}}\mathcal{G}_\theta(\mathbf{u}_S)^\top$.
In particular, the first term $\nabla_{\mathbf{u}}\mathcal{G}_\theta(\mathbf{u}_S)^\top[\cdot]$ captures the \emph{indirect} effect of $\mathbf{u}_S$ on the objective mediated by the state, while the second bracket collects \emph{direct} control contributions.
The penalty induces an implicit multiplier $\boldsymbol{\nu}_\mu=2\mu\mathcal{R}_\theta(\mathbf{u}_S)$, and, consequently, the residual term influences the gradient in the same manner as a Lagrange multiplier influences the adjoint/KKT system.
\begin{remark}
Proposition~\ref{prop:implicit_adjoint} is a chain-rule factorization of the \emph{reduced} gradient. 
The bracketed terms coincide with the adjoint equation residual \eqref{eq:adjoint} and Lagrangian gradient \eqref{eq:adjoint_gradient} evaluated at 
$(\mathbf{y},\mathbf{u},\boldsymbol{\nu})=(\mathcal{G}_\theta(\mathbf{u}_S),\mathbf{u}_S,2\mu\mathcal{R}_\theta(\mathbf{u}_S))$, 
but they do not need to vanish individually at a reduced stationary point since $\mathbf{y}$ is not an independent variable.
\end{remark}

When $\mu=0$, the implicit multiplier vanishes, $\boldsymbol{\nu}_\mu\equiv \mathbf{0}$, and the reduced gradient in
\eqref{eq:implicit_adjoint_grad} reduces to the tracking gradient $\nabla_\mathbf{u}\bar{J}_{NO}$. Achieving a feasible minimum for pure tracking via gradient descent is notoriously ill-posed in practice. In the context of a surrogate NO, this ill-posedness manifests as the optimizer aggressively minimizes $\bar{J}_{NO}$ by exploiting network approximation errors, finding non-physical inputs that the NO erroneously maps to the target, even when additional Tikhonov regularization is used. This numerical reality strictly necessitates the residual penalty $\mu\|\mathcal{R}_\theta\|^2$, not just to approximate the adjoint, but to actively block these non-physical optimization directions, which we formalize next as \textit{cheating directions}.

\subsection{Cheating Directions}\label{section3.4}

Optimizing only on the cost function, i.e., tracking, allows the optimizer to exploit errors in the NO to reduce the tracking cost by exploring directions that deviate from the PDE mapping. We call these directions \textit{cheating directions}, and define them formally as follows.

\begin{definition}[Cheating direction]\label{def:cheating_direction}
Let $\bar{J}_{\mathrm{NO}}$ be the tracking objective and
$R(\mathbf{u}_S):=\|\mathcal{R}_\theta(\mathbf{u}_S)\|_2^2$ the residual penalty.
At a point $\mathbf{u}_S$, a direction $\mathbf{d}$ is called a \emph{cheating direction} if
\[
\nabla_{\mathbf{u}}\bar{J}_{\mathrm{NO}}(\mathbf{u}_S)^\top \mathbf{d} < 0
\quad\text{and}\quad
\nabla_{\mathbf{u}}R(\mathbf{u}_S)^\top \mathbf{d}> 0.
\]
\end{definition}

We show that these cheating directions exist in the following proposition for the reduced problem \eqref{eq:pde_simple}.
\begin{proposition}[Existence of cheating directions]\label{prop:cheating_directions}
Let $R(\mathbf{u}_S)=\|\mathcal{R}_\theta(\mathbf{u}_S)\|^2$ and define $\nabla_\mathbf{u}\bar{J}_{NO}(\mathbf{u}_S):=\mathbf{g}$ and $\nabla_\mathbf{u}R(\mathbf{u}_S):=\mathbf{r}$ as the gradients of the cost function and the residual penalty term, respectively. If $\mathbf{g}\neq0$, $\mathbf{r}\neq0$, and $\mathbf{g}$ and $\mathbf{r}$ are not collinear, i.e. $\mathbf{g}\neq c\mathbf{r}$ for some constant $c\neq0$, then there exists a cheating direction $\mathbf{d}$ such that
\begin{equation}\label{eq:cheating_direction}
\mathbf{d}^\top \mathbf{g}<0 \quad \text{and}\quad \mathbf{d}^\top \mathbf{r}>0.
\end{equation}
Consequently, for sufficiently small $\eta>0$, it follows that
\begin{equation}\label{eq:cheating_direction_step}
\bar{J}_{NO}(\mathbf{u}_S+\eta \mathbf{d})<\bar{J}_{NO}(\mathbf{u}_S) \quad \text{and}\quad R(\mathbf{u}_S+\eta \mathbf{d})>R(\mathbf{u}_S).
\end{equation}
\end{proposition}

\begin{proof}
We consider three cases based on the sign of 
$\mathbf{g}^\top\mathbf{r}$.

\medskip\noindent
\textbf{Case 1: $\mathbf{g}^\top \mathbf{r}>0$.}
Consider the family of directions
\begin{equation}\label{eq:family_pos}
\mathbf{d}=-\mathbf{g}+\alpha \mathbf{r}, \quad \alpha\in \mathbb{R}.
\end{equation}
Then
\begin{equation}
\mathbf{g}^\top \mathbf{d}=\mathbf{g}^\top(-\mathbf{g}+\alpha \mathbf{r})=-\|\mathbf{g}\|^2+\alpha \mathbf{g}^\top \mathbf{r},
\end{equation}
and similarly
\begin{equation}
\mathbf{r}^\top \mathbf{d}=\mathbf{r}^\top(-\mathbf{g}+\alpha \mathbf{r})=-\mathbf{r}^\top \mathbf{g}+\alpha\|\mathbf{r}\|^2
=-(\mathbf{g}^\top \mathbf{r})+\alpha\|\mathbf{r}\|^2.
\end{equation}
Thus $\mathbf{g}^\top \mathbf{d}<0$ holds whenever $\alpha<\|\mathbf{g}\|^2/(\mathbf{g}^\top \mathbf{r})$, and $\mathbf{r}^\top \mathbf{d}>0$ holds whenever
$\alpha>(\mathbf{g}^\top \mathbf{r})/\|\mathbf{r}\|^2$. Such an $\alpha$ exists if and only if
\begin{equation}
\frac{\mathbf{g}^\top \mathbf{r}}{\|\mathbf{r}\|^2}<\frac{\|\mathbf{g}\|^2}{\mathbf{g}^\top \mathbf{r}}
\iff (\mathbf{g}^\top \mathbf{r})^2 <\|\mathbf{g}\|^2\|\mathbf{r}\|^2,
\end{equation}
which holds by the strict Cauchy--Schwarz inequality since $\mathbf{g}$ and $\mathbf{r}$ are not collinear.
Therefore there exists $\alpha$ such that $\mathbf{g}^\top \mathbf{d}<0$ and $\mathbf{r}^\top \mathbf{d}>0$.

\medskip\noindent
\textbf{Case 2: $\mathbf{g}^\top \mathbf{r} < 0$.}
Take $\mathbf{d} = -\mathbf{g}$. Then
\begin{equation}
\mathbf{g}^\top\mathbf{d} = -\|\mathbf{g}\|^2 < 0,
\end{equation}
and
\begin{equation}
\mathbf{r}^\top\mathbf{d} = -\mathbf{g}^\top\mathbf{r} > 0.
\end{equation}
That is, the steepest descent direction for the tracking cost is itself
a cheating direction.

\medskip\noindent
\textbf{Case 3: $\mathbf{g}^\top \mathbf{r} = 0$.}
Take $\mathbf{d} = -\mathbf{g} + \alpha\mathbf{r}$ with any $\alpha > 0$. Then
\begin{equation}
\mathbf{g}^\top\mathbf{d} = -\|\mathbf{g}\|^2 < 0,
\end{equation}
and
\begin{equation}
\mathbf{r}^\top\mathbf{d} = \alpha\|\mathbf{r}\|^2 > 0.
\end{equation}

\medskip\noindent
Finally, since $\bar{J}_{\mathrm{NO}}$ and $R$ are differentiable, for sufficiently small $\eta>0$ we have
\[
\bar{J}_{\mathrm{NO}}(\mathbf{u}_S+\eta \mathbf{d})=\bar{J}_{\mathrm{NO}}(\mathbf{u}_S)+\eta\, \mathbf{g}^\top \mathbf{d}+o(\eta),
\qquad
R(\mathbf{u}_S+\eta \mathbf{d})=R(\mathbf{u}_S)+\eta\, \mathbf{r}^\top \mathbf{d}+o(\eta),
\]
so $\mathbf{g}^\top \mathbf{d}<0$ and $\mathbf{r}^\top \mathbf{d}>0$ imply
$\bar{J}_{\mathrm{NO}}(\mathbf{u}_S+\eta \mathbf{d})<\bar{J}_{\mathrm{NO}}(\mathbf{u}_S)$ and
$R(\mathbf{u}_S+\eta \mathbf{d})>R(\mathbf{u}_S)$ for all sufficiently small $\eta>0$.
\end{proof}

Proposition \ref{prop:cheating_directions} demonstrates that there always exist directions where the tracking loss can be improved at the expense of the residual during the optimization iterations, if the gradients are not zero nor collinear.  Since the gradients $\mathbf{g}$ and $\mathbf{r}$ are multidimensional gradients of different functions, they are not expected to be aligned, which is what we consistently observe in our computational experiments. 

\begin{remark}
    Cheating directions are analogous to adversarial perturbations in machine learning \cite{goodfellow2015adversarial}, where small input modifications exploit model approximation errors to produce outputs that minimize a loss while deviating from the true input–output mapping. Here, the optimizer plays the role of an adversary, finding controls that reduce the tracking cost by exploiting regions where the NO's approximation of the PDE solution operator is inaccurate, rather than by following the true PDE dynamics.
\end{remark} 

Including a residual penalty in the cost function prevents the optimizer from exploiting cheating directions. The residual penalty can thus be interpreted as an adversarial robustness mechanism for the surrogate optimization. For a sufficiently large penalty factor $\mu$, the optimizer does not exploit the cheating direction. This result is summarized in the following corollary.

\begin{corollary}[Suppressing the cheating directions]\label{cor:suppress_cheating}
Let $\bar{J}_\mu:=\bar{J}_{NO}(\mathbf{u}_S)+\mu R(\mathbf{u}_S)$ and $\mathbf{g}$ and $\mathbf{r}$ be defined as in Proposition \ref{prop:cheating_directions}. Suppose $\mathbf{d}$ is a cheating direction at $\mathbf{u}_S$, i.e.,
\begin{equation}
\mathbf{d}^\top \mathbf{g}<0 \quad \text{and}\quad \mathbf{d}^\top \mathbf{r}>0.
\end{equation}
Then there exists a threshold $\mu^*$ as
\begin{equation}
\mu^*:=\frac{-\mathbf{g}^\top \mathbf{d}}{\mathbf{r}^\top \mathbf{d}}>0
\end{equation}
such that for all $\mu>\mu^*$,
\begin{equation}
(\nabla_\mathbf{u}\bar{J}_\mu)^\top \mathbf{d}=\mathbf{g}^\top \mathbf{d} + \mu \mathbf{r}^\top \mathbf{d}>0.
\end{equation}
\end{corollary}
\begin{proof}
The proof follows directly by solving for $\mu$ from $(\nabla_\mathbf{u}\bar{J}_\mu)^{\top}\mathbf{d}>0$.
\end{proof}

The consequence of Corollary \ref{cor:suppress_cheating} is that the directional derivative of the penalized objective along a cheating direction becomes positive once $\mu$ is large enough, and the cheating direction $\mathbf{d}$ is no longer a descent direction for the penalized objective, preventing the optimizer from exploiting that direction. In practice, for a sufficiently large $\mu$, the penalty factor can be kept constant during the iterative optimization process, as we demonstrate in our experiments in section \ref{section_5}.

In our experiments, we observe that cheating directions change frequently. Consequently, the penalty term in \eqref{eq:pde_penalized_reduced} also acts as a regularizer, typically enforcing smoothness and suppressing high-frequency content in the controls in case the underlying PDE admits smooth solutions, and the control is additive to the PDE. That is, adding the residual penalty suppresses the constantly changing cheating directions, leading to observable improved stability of the optimization and markedly lower high-frequency content in the optimized controls. 

\begin{remark}\label{rem:filtering}
For many PDEs with elliptic or diffusive character (e.g., Poisson, viscous Burgers, steady Navier--Stokes with viscosity), the forward map from forcing/control to state suppresses high-frequency components. Thus, a control dominated by high frequencies typically produces a comparatively small change in the state, while it can strongly affect the residual when evaluated through the operator. Consequently, enforcing a small residual through $\mu R(u)$ discourages the optimizer from adding high-frequency content to $u$ that is not supported by the PDE dynamics, behaving like an implicit low-pass regularizer.
\end{remark}
\begin{remark}\label{remark:smoothness}
For example, consider the 1D Poisson's problem $\Delta y(x)=u(x)$. Taking the Fourier transform on both sides yields $\xi^2\widehat{y}(\xi)=\widehat{u}(\xi)$, which can be rearranged to $\widehat{y}(\xi)=\frac{\widehat{u}(\xi)}{\xi^2}$. Thus, $|\widehat{y}(\xi)|$ decays at rate $1/\xi^2$, thereby effectively filtering high frequency components in the state $y$. Furthermore, the PDE residual $\Delta y(x)-u(x)$ contains components of $\xi^2$ in the Fourier space, and therefore implicitly places a high cost on high-frequency components, which explains our empirical observations and suppression of the oscillatory cheating directions.
\end{remark}

\begin{remark}\label{rem:regularization}
The regularizing term $\|\mathbf{u}_S\|_{\mathcal{U}_h}^2$ may help in avoiding high frequency (erroneous) solutions and the cheating directions, if it enforces smoothness, for example, with an $H^1$ or $H^2$ seminorm. However, the regularizer does not enforce the differential equation constraints on the optimizer, and thus $\mathcal{R}_\theta$ can be large. If the PDE is not strongly smoothing (e.g., predominantly hyperbolic or weakly dissipative), adding $H^1/H^2$ regularization can complement the residual penalty, as we show in \ref{appendix_d}.
\end{remark}

\section{EXPERIMENTAL SETUP} \label{section_4}
We conduct experiments of physics-informed DeepONet models for control problems of tracking type, showcasing our proposed method's re-usability, how it extends the original scope of NOs, and its architectural simplicity as well as performance as an effective surrogate model. For each PDE, we define two continuous tracking objectives $J_{i}(\mathbf{u})$ for $i=1,2$. In all experiments, this is approximated by quadrature on a uniform $(32,32)$ sensor grid $S$, and represent the control $\mathbf{u}$ by its values on $S$, that is $\mathbf{u}_S:= [\mathbf{u}(\mathbf{s}_1),\dots,\mathbf{u}(\mathbf{s}_N)]^\top$, and optimized via the NO-reduced penalized objective $\bar{J}_\mu(\mathbf{u}_S)=\bar{J}_{NO}(\mathbf{u}_S)+\mu\|\mathcal{R}_\theta(\mathbf{u}_S)\|_2^2$, as described in \eqref{eq:pde_tracking_no_penalized_reduced}. We further consider the problems to be reachable, meaning that the control is capable of reaching the desired tracked state. Reachability allows us to create ground-truth states with a known control, and thus, to perform a numerical a posteriori error analysis.

\subsection{Experimental Problems}\label{section_4.1}

We selected problems that are common PDEs used in both test cases in the literature and practical applications. For each problem, we list two different cost functions and solve the arising tracking problems using our method, demonstrating that the trained physics-informed DeepONet can be used directly as a surrogate model. We study the following PDE systems. 

\begin{enumerate}
  \item \textit{Scalar Elliptic Control: Poisson's Problem}
    \begin{equation}\label{eq:poisson}
    \begin{aligned}
      \Delta y(\mathbf{x}) &= -u(\mathbf{x}), \qquad \mathbf{x}\in\Omega,\\
      y(\mathbf{x}) &= 0, \qquad \mathbf{x}\in\partial\Omega,\\
      \Omega &= [0,1]\times[0,1].
    \end{aligned}
    \end{equation}
This problem is a test case from \cite{Hwang2022} and reflects targeting a heat profile with a given source. A slightly modified version can also be found in \cite{Yong2024Deep}. The first cost function follows a state that was generated with $u_{ref,1}(x_1,x_2)=\sin(\pi x_1)\sin(\pi x_2)$, which is a scaled version from \cite{Hwang2022}. The second tracking object is a solution obtained by a control generated by a Gaussian random field (GRF) with a lengthscale of $l=1.0$.
\begin{equation}\label{cost:poisson}
    J_{\mathrm{Pois},i}(u) = \int_{\Omega} \big(y(\mathbf{x}) - y_{d,i}(\mathbf{x})\big)^2\,d\mathbf{x} + \lambda\|u\|^2_{\mathcal{U}},  \quad i=1,2.
  \end{equation}

  \item \textit{Nonlinear Transport Control: Viscous Burgers' Equation}
    \begin{equation}\label{eq:burgers}
    \begin{aligned}
    \frac{\partial y(x,t)}{\partial t} + y\,\frac{\partial y(x,t)}{\partial x}
    &= 0.01\,\frac{\partial^2 y(x,t)}{\partial x^2} + u(x,t),\\
    y(x,0) &= 0,\quad y|_{\partial\Omega}=0,\\
    (x,t) \in \Omega&\times[0,T],\ \Omega=[0,1].
    \end{aligned}
    \end{equation}
We choose this problem in order to demonstrate that our method works for nonlinear time-dependent PDEs. A variant of this problem was given in \cite{Hwang2022}. We choose two cost functions of tracking type, where the tracked references are solutions of controls generated by a GRF with lengthscale of $l=1.0$ for $J_{\mathrm{Burg},1}$ and $l=0.5$ for $J_{\mathrm{Burg},2}$.  
\begin{equation}\label{cost:burgers}
J_{\mathrm{Burg},i}(u)
= \int_0^T\!\!\int_{\Omega} \big(y(x,t) - y_{d,i}(x,t)\big)^2\,dx\,dt + \lambda \|u\|^2_{\mathcal{U}}, \quad \ i=1,2.\\
\end{equation}

\item \textit{Flow Control: The Stokes Equation}
\begin{equation}\label{eq:stokes}
\begin{aligned}
0.1\Delta \mathbf{v}(\mathbf{x}) + \nabla p(\mathbf{x})
&= \mathbf{u}(\mathbf{x}), 
\qquad \mathbf{x}\in\Omega,\\
\nabla \cdot \mathbf{v}(\mathbf{x})
&= 0,
\qquad \mathbf{x}\in\Omega,\\
\mathbf{v}(\mathbf{x})
&= \bigl(2x_2(1-x_2),\,0\bigr)^\top,
\qquad \mathbf{x}\in\Gamma_{\mathrm{in}},\\
\mathbf{v}(\mathbf{x})
&= \mathbf{0},
\qquad \mathbf{x}\in\Gamma_{\mathrm{w}},\\
\frac{\partial \mathbf{v}}{\partial n_1}(\mathbf{x})
&= \mathbf{0},
\qquad \mathbf{x}\in\Gamma_{\mathrm{out}},\\
\int_{\Omega} p(\mathbf{x})\,d\mathbf{x}
&= 0. 
\end{aligned}
\end{equation}

    where the domain is $\Omega=[0,1]\times[0,1]$. The boundary is partitioned into inlet $\Gamma_{\mathrm{in}}=\{0\}\times[0,1]$, walls $\Gamma_{\mathrm{w}}=[0,1]\times\{0,1\}$, and outlet $\Gamma_{\mathrm{out}}=\{1\}\times[0,1]$. The inlet profile is a fixed parabola with amplitude $0.5$, the walls enforce no-slip, and the outlet uses a homogeneous Neumann condition on the velocity. We include $\int_{\Omega}p(\mathbf{x})\,d\mathbf{x}=0$ to fix the pressure up to a constant. We choose this problem to demonstrate that our method also works for vector controls. An alternative version of this problem can be found in \cite{Demo2023}. The objective is to recover the control (disturbance) field $\mathbf{u}$ that reproduces a desired reference velocity $\mathbf{v}_{d,i}$. The reference controls were drawn from a GRF with lengthscale $l=1.0$ for $J_{\mathrm{Stokes},1}$ and $l=0.5$ for $J_{\mathrm{Stokes},2}$, and the corresponding tracked reference velocities were obtained by solving the forward problem with our reference solver. We consider a tracking cost with a regularizer:
    \begin{equation}\label{cost:stokes}
    J_{\mathrm{Stokes},i}(\mathbf{u})
      = \int_{\Omega} \|\mathbf{v}(\mathbf{x}) - \mathbf{v}_{d,i}(\mathbf{x})\|^2\,d\mathbf{x} + \lambda\|\mathbf{u}\|^2_{\mathcal{U}}, \quad i=1,2.\\
    \end{equation}

\end{enumerate}

\subsection{Network Architectures and Training.}
We employed the DeepONet architecture, using a modified version with residual connections in the fully connected networks as described by Wang et al. \cite{wang2021understanding}. The architecture details, hyperparameters, and training loss curves are given in \ref{appendix_a}. Training was purely physics-informed, i.e., without a data loss, with the loss consisting of the differential residual \eqref{residual}, boundary residual \eqref{residual_boundary}, and initial condition residual \eqref{residual_initial}.

\subsection{Regularization}
In order to demonstrate the existence of cheating directions as in Proposition \ref{prop:cheating_directions}, we solve the scalar elliptic control and the nonlinear transport control with a neural operator approach (as in \eqref{eq:pde_tracking_no_reduced}) for two different regularizers without the residual penalty ($\mu=0)$ and compare against the residual penalty with no additional regularizer ($\lambda=0$). We use the squared $L_2$-norm and the squared $H^1$-seminorm as regularizers. By using the common $L_2$ regularizer, we show that the cheating directions are easily exploitable by an optimizer. Further, as noted in Remark \ref{remark:smoothness}, adding a smoothing regularizer, such as $H^1$ may block to some extent cheating directions, but it is not sufficient for achieving a feasible and accurate solution, which can be observed in our experiments.

In the experiments, $\|\cdot\|_{\mathcal U}^2$ denotes either
$\|\mathbf{u}\|_{L_2}^2$ or $|\mathbf{u}|_{H^1}^2:=\|\nabla \mathbf{u}\|_{L_2}^2$,
with corresponding discrete realizations on the sensor grid. In discrete form, we use
\begin{equation}
\|\mathbf{u}_h\|_{\mathcal{U}_h}^2 := \sum_{j=1}^{m} \|\mathbf{u}(\mathbf{x}_j,t_j)\|_2^2\,w_j,
\quad
|\mathbf{u}_h|_{H^1}^2 := \sum_{j=1}^{m} \|\nabla_h \mathbf{u}(\mathbf{x}_j,t_j)\|_2^2\,w_j,
\end{equation}
where $w_j$ are quadrature weights and $\nabla_h$ denotes a discrete first derivative operator, and is applied componentwise when $\mathbf{u}$ is vector-valued. For the general quadrature, we use a simple Riemann integral. For evaluating the discrete gradients, we use the forward difference for $H^1$.

\subsection{Baseline Optimization Parameters}
We use Adam with decoupled weight decay optimizer (AdamW) provided by Optax \cite{deepmind2020jax}. We use the default hyperparameters $b_1=0.9$ and $b_2=0.999$ for the exponential moving averages of the first and second moments, a numerical stabilizer $\varepsilon=10^{-8}$, and $\varepsilon_{\mathrm{root}}=0.0$. We disable decoupled weight decay to avoid introducing additional $L_2$ regularization. Moreover, we use no parameter masking, and no Nesterov-style momentum. The initial update step size (learning rate) is set as $\gamma=0.1$, and we use a decay rate of $0.5$ at each 200 steps. This prevents oscillations after the optimizer has settled on a solution.  

We solve the different problems using the NO-reduced formulation \eqref{eq:pde_tracking_no_penalized_reduced} and evaluate two stabilizing mechanisms: 
(i) \emph{residual penalization only} ($\mu>0$, $\lambda=0$), and 
(ii) \emph{regularization only} ($\mu=0$, $\lambda>0$), where $\mu$ and $\lambda$ are defined in \eqref{eq:pde_penalized_reduced}. 
For each problem, the nonzero parameter (either $\mu$ or $\lambda$) was selected by a small grid search to minimize the control error while keeping all other settings fixed. 
The optimizer's number of iterations $n_{\mathrm{iter}}$ is identical across all runs. We summarize the baseline parameters in Table \ref{tab:parameters_new}

\begin{table}[h!]
\centering
\caption{Optimization settings used in the experiments (applied for both targets $i=1,2$). We compare residual-penalty-only runs ($\mu>0,\lambda=0$) and regularization-only runs ($\mu=0,\lambda>0$) for the scalar elliptic control and the nonlinear transport control.}
\label{tab:parameters_new}
\scriptsize
\begin{tabular}{lccc|ccc}
\toprule
\textbf{Problem} & $n_{\mathrm{iter}}$ & $\gamma_0$ & decay & $\mu$ (penalty-only) & $L_2: \lambda$ (reg-only) &  $H^1:\lambda$ (reg-only)\\
\midrule
$J_{\mathrm{Pois},i}(u)$   & 2000 & 0.1 & $\times 0.5/200$ & 0.01 & $1\times 10^{-4}$ &$1\times 10^{-4}$ \\
$J_{\mathrm{Burg},i}(u)$   & 2000 & 0.1 & $\times 0.5/200$ & 0.01 & $1\times 10^{-2}$ & $1\times 10^{-2}$ \\
$J_{\mathrm{Stokes},i}(\mathbf{u})$ & 2000 & 0.1 & $\times 0.5/200$ & 0.1  & - \\
\bottomrule
\end{tabular}
\end{table}

\subsection{Reference Solutions}
As reference solutions, we employ the adjoint method for the scalar elliptic and the nonlinear transport control problems. We implement the forward PDE-solvers and adjoint PDE-solvers in JAX \cite{jax2018github}. The details are given in \ref{appendix_b}. We use the same grid size, collocation points (sensor locations), and optimizer with the same update step size $\gamma$ and scheduler as in our NO approach. This allows us to make a fair comparison of solution times, time per iteration, and accuracy. 

For the flow control problem, we do not implement an adjoint-based solver. The tracking objective \eqref{cost:stokes} considers only the velocity field, leaving the pressure undetermined. Consequently, the control-to-velocity map is not injective: multiple controls can generate the same target velocity with different pressure fields, rendering a pointwise posterior-error analysis of the control ill-defined. Furthermore, the adjoint system for the flow control involves a saddle-point structure that is substantially more complex to implement than the scalar adjoint systems arising in the scalar elliptic and the nonlinear transport control problem. A rigorous adjoint-based comparison for Stokes is left to future work. Instead, we verify the feasibility of the optimized control by solving the forward Stokes equations with an independent finite difference solver and comparing the resulting velocity fields.

We compare the adjoint method with two regularizers, $L_2$ and $H^1$, against the penalty method in terms of accuracy and computational speed. We select the best parameters that gave the smallest mean square error (MSE) for the control when compared against the true solution. For obtaining the best regularization parameter, we did a parameter sweep for each problem with regularization values $\lambda=\{10^{-10},10^{-9},\dots,10^{-2}\}$. The results of the sweep are shown in \ref{appendix_b_regularization} and the best values are summarized in Table \ref{tab:ref_best_lambda_summary}.

\begin{table}[!h]
\centering
\caption{Best adjoint method regularization parameters $\lambda$ for $L_2$ and $H^1$ regularization, selected by the smallest control mean squared error.}
\label{tab:ref_best_lambda_summary}
\renewcommand{\arraystretch}{1.15}
\setlength{\tabcolsep}{5pt}

\begin{tabular}{c cccc}
\toprule
Method & $J_{Pois,1}$ & $J_{Pois,2}$ & $J_{Burg,1}$ & $J_{Burg,2}$ \\
\midrule
$L_2$ & $\lambda=10^{-5}$ & $\lambda=10^{-8}$ & $\lambda=10^{-3}$ & $\lambda=10^{-4}$ \\
$H^1$ & $\lambda=10^{-9}$ & $\lambda=10^{-9}$ & $\lambda=10^{-5}$ & $\lambda=10^{-9}$ \\
\bottomrule
\end{tabular}
\end{table}

All methods (including the neural operator) are executed on the GPU. We ran the problems on a system with 12th Gen Intel(R) Core(TM) i5-12600K, 32GB RAM and an NVIDIA RTX(TM) A2000 12GB.

\section{RESULTS}\label{section_5}

We present our results visually in the following subsections, followed by a brief sensitivity study of the penalty and initial step-size parameters $\mu$ and $\gamma$ for the scalar elliptic control and the nonlinear transport control problems. For the same problems, we summarize solution time, time per iteration, and mean square error (MSE) of the control error against the reference solution in the final section.

\subsection{Scalar Elliptic Control: The Poisson Equation}

\begin{figure}[!t]
    \centering

    \begin{subfigure}{\linewidth}
        \centering
        \includegraphics[width=\linewidth]{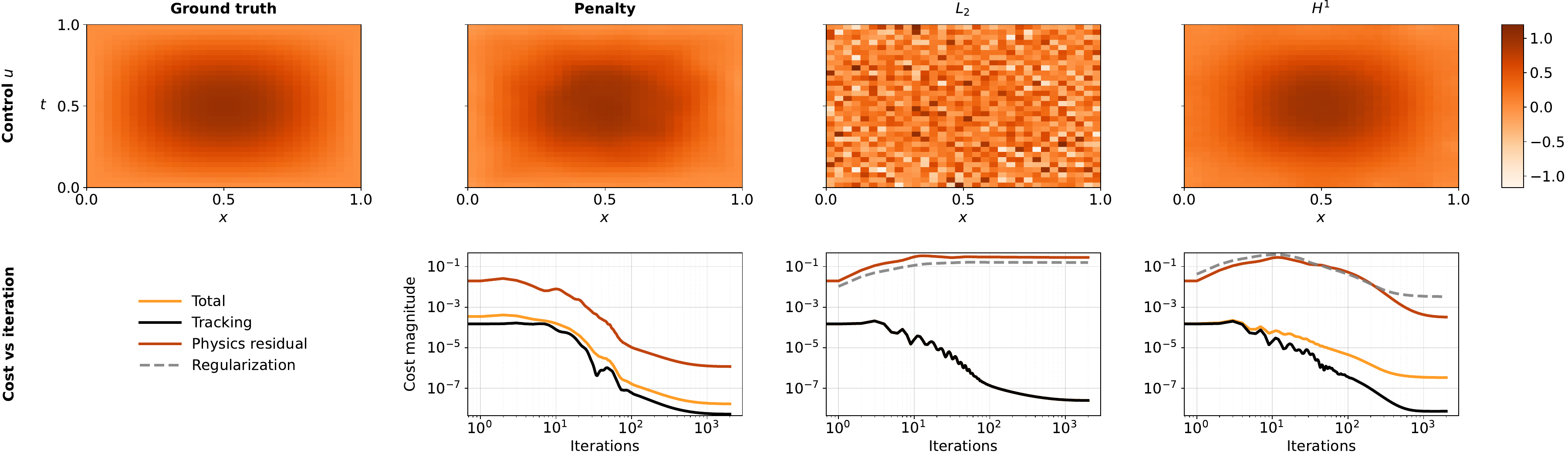}
        \caption{Objective $J_{\mathrm{Pois},1}$.}
        \label{fig:poisson_1_regularization}
    \end{subfigure}

    \vspace{1pt}

    \begin{subfigure}{\linewidth}
        \centering
        \includegraphics[width=\linewidth]{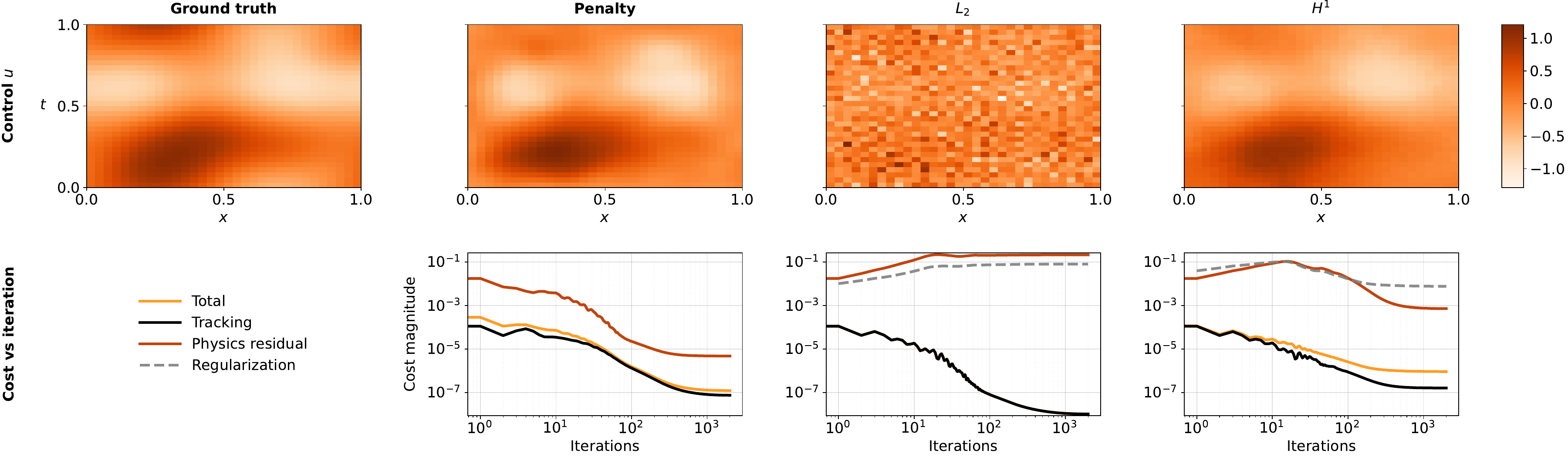}
        \caption{Objective $J_{\mathrm{Pois},2}$.}
        \label{fig:poisson_2_regularization}
    \end{subfigure}

    \caption{Scalar elliptic control problem \eqref{eq:poisson} with costs \eqref{cost:poisson}: comparison of three NO-approaches against the ground truth. Columns (left to right) show: ground truth (reference), \emph{penalty-only} method ($\lambda=0$), \emph{$L_2$-regularization only} ($\mu=0$), and \emph{$H^1$-regularization only} ($\mu=0$). The top row shows the final control $u$ and the bottom row the optimization cost versus iteration.}
    \label{fig:poisson_regularization}
\end{figure}

\begin{figure}[!t]
\centering
\begin{subfigure}{\linewidth}
  \centering
  \includegraphics[width=0.9\linewidth]{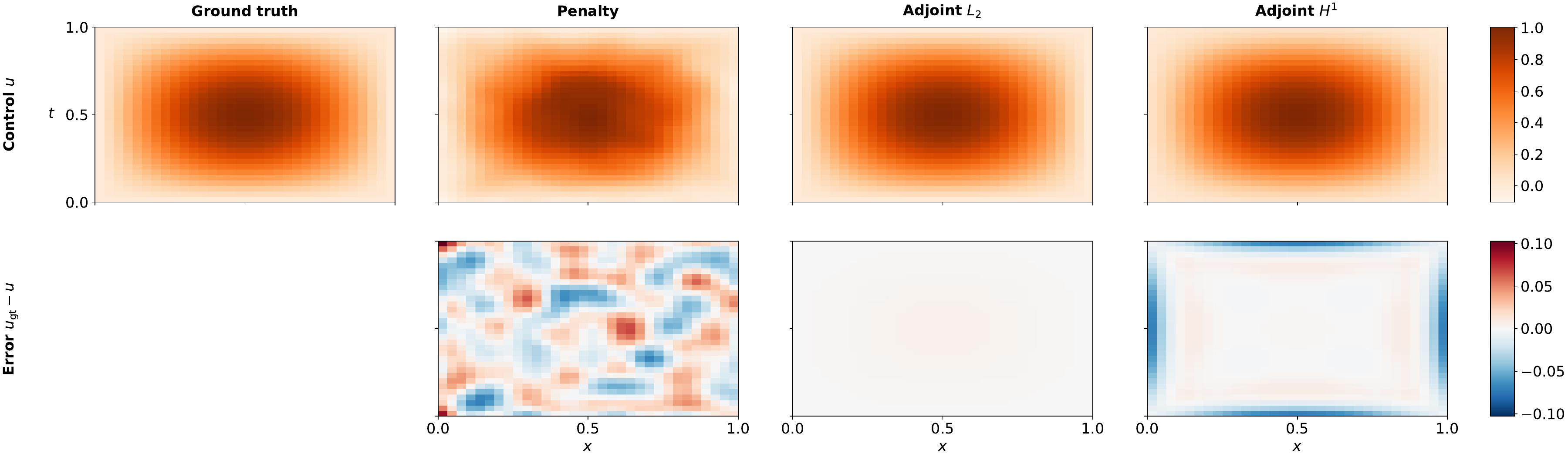}
    \caption{Objective $J_{Pois,1}$.}
  \label{fig:poisson_1_regularization_reference}
\end{subfigure}

\vspace{1pt}

\begin{subfigure}{\linewidth}
  \centering
  \includegraphics[width=0.9\linewidth]{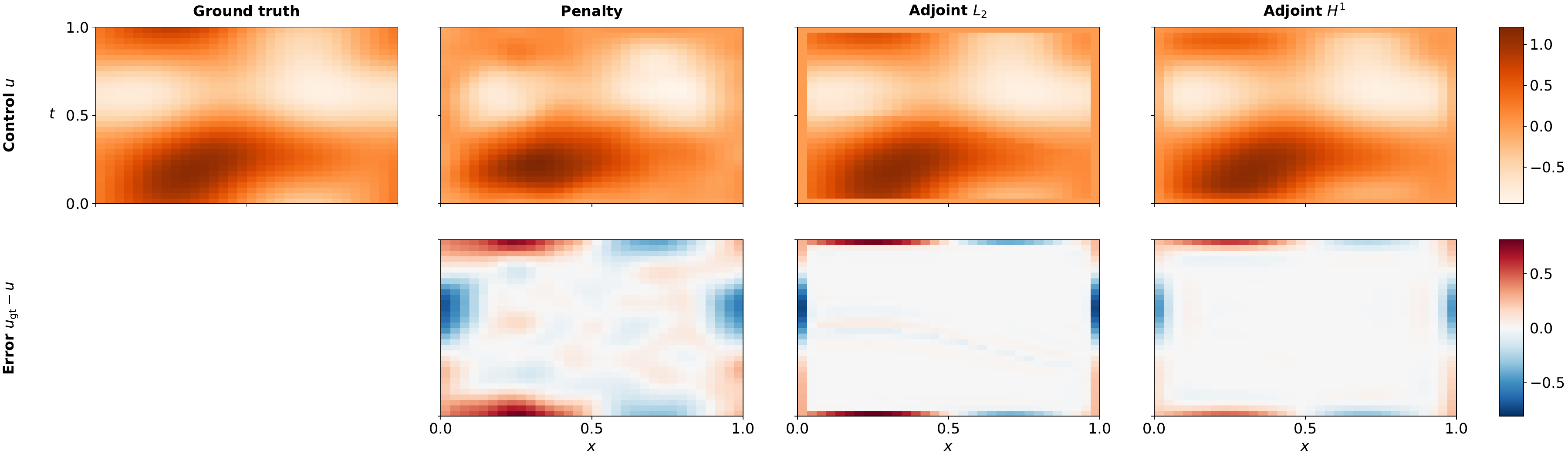}
        \caption{Objective $J_{Pois,2}$.}
  \label{fig:poisson_2_regularization_reference}
\end{subfigure}

\caption{Scalar elliptic control problem \eqref{eq:poisson} with costs \eqref{cost:poisson}: comparison of three approaches against the ground truth. Columns (left to right) show: ground truth (reference), \emph{penalty-only} method ($\lambda=0$), adjoint method with \emph{$L_2$-regularization}, and adjoint method with \emph{$H^1$-regularization}. The top row shows the final control $u$ and the bottom row the control error $u-u_{ref}$.}
\label{fig:poisson_regularization_reference}
\end{figure}

We consider the tracking problem \eqref{cost:poisson} with the Poisson's equation \eqref{eq:poisson} as the underlying constraint. We first introduce the effects of the cheating directions by using only $L_2$ regularization and without penalty. The effects are clearly visible in Figures \ref{fig:poisson_1_regularization} and \ref{fig:poisson_2_regularization} for both problems $J_{Pois,1}$ and $J_{Pois,2}$, where the optimizer finds a noisy control that produces a state that matches the tracking target by exploiting the cheating directions. The obtained solution has a high physics residual in the magnitude of $\mathcal{O}(1)$ and does not satisfy the PDE constraints; thus, it is not a feasible solution. We tried to alter the regularization parameter within a range of $\lambda\in[10^{-8},1]$, but found no difference in the quality of the solution. However, a smoothing regularization $H^1$ produces a control that is smooth but has a doubly higher physics residual than the penalized method, thus violating the constraint. A summary of the results is given in Table \ref{tab:summary_errors_poisson}. The penalty method shows the best overall performance regarding tracking error and the feasibility, i.e., the physics residual.

We further compare our method against the adjoint method, where we solve the adjoint equations with a finite difference scheme. As before, we employ  $L_2$ and $H^1$ as regularization. Figures \ref{fig:poisson_1_regularization_reference} and \ref{fig:poisson_2_regularization_reference} shows the optimized controls $u_S(\mathbf{x})$, and the error $u_S(\mathbf{x})-u_{ref}(\mathbf{x})$ where $u_{ref,i}(\mathbf{x})$ generated the tracked state $y_{d,i}(\mathbf{x})$. From the Figure \ref{fig:poisson_2_regularization_reference}, we observe that the adjoint method produces accurate results on the interior, but does not find a proper solution at the boundaries. The reason is that the Dirichlet boundary conditions do not enter the PDE solution process as variables (degrees of freedom) and, thus, do not influence the solution. Our penalized NO also shows deficient performance on the boundary due to the same reasons, i.e., control values on the boundary do not influence the solution. For the case $J_{Pois,1}$, as shown in Figure \ref{fig:poisson_1_regularization_reference}, the boundary error effect is not present for the $L_2$ and $H^1$ regularization as the control takes a value of 0 on the boundary and the initial starting guess for optimization was zero, i.e., $u_{S,0}(\mathbf{x})=0$. Further, our method shows inferior performance on the interior than the adjoint method. This is naturally due to Poisson's problem and its adjoint equation being a linear PDE. Therefore, the adjoint method can solve it more accurately than a surrogate model.

\begin{table}[t]
\centering
\scriptsize
\caption{Summary of errors for scalar elliptic control problem for different penalties, $L_2$, and $H^1$ regularization.}
\label{tab:summary_errors_poisson}
\renewcommand{\arraystretch}{1.15}
\setlength{\tabcolsep}{3pt}

\begin{tabular}{c l ccc c ccc}
\toprule
& & \multicolumn{3}{c}{Penalty method} & \multicolumn{1}{c}{$L_2$ regularization} & \multicolumn{3}{c}{$H^1$ regularization} \\
\cmidrule(lr){3-5}\cmidrule(lr){6-6}\cmidrule(lr){7-9}
Problem & Metric & $\mu=10^{-1}$ & $\mu=10^{-2}$ & $\mu=10^{-3}$ & $\lambda=10^{-4}$ & $\lambda=10^{-3}$ & $\lambda=10^{-4}$ & $\lambda=10^{-5}$ \\
\midrule
\multirow{3}{*}{$J_{Pois,1}$} & Tracking error $\mathcal{J}_{\text{track}}$ & $2.97\times 10^{-8}$ & $5.36\times 10^{-9}$ & $1.01\times 10^{-8}$ & $2.53\times 10^{-8}$ & $4.68\times 10^{-7}$ & $7.47\times 10^{-9}$ & $1.57\times 10^{-8}$ \\
 & Physics residual $\mathcal{R}_\theta$ & $1.2\times 10^{-6}$ & $1.19\times 10^{-6}$ & $1.02\times 10^{-5}$ & $2.81\times 10^{-1}$ & $1.13\times 10^{-4}$ & $3.24\times 10^{-4}$ & $6.61\times 10^{-2}$ \\
 & Control error (MSE) & $1.21\times 10^{-3}$ & $6.29\times 10^{-4}$ & $2.61\times 10^{-3}$ & $2.72\times 10^{-1}$ & $4.53\times 10^{-2}$ & $3.95\times 10^{-3}$ & $6.63\times 10^{-2}$ \\
\midrule
\multirow{3}{*}{$J_{Pois,2}$} & Tracking error $\mathcal{J}_{\text{track}}$ & $1.19\times 10^{-6}$ & $7.59\times 10^{-8}$ & $3.88\times 10^{-8}$ & $1.04\times 10^{-8}$ & $1.6\times 10^{-6}$ & $1.63\times 10^{-7}$ & $4.73\times 10^{-8}$ \\
 & Physics residual $\mathcal{R}_\theta$ & $2.95\times 10^{-6}$ & $4.82\times 10^{-6}$ & $1.16\times 10^{-5}$ & $2.15\times 10^{-1}$ & $1.95\times 10^{-4}$ & $7.24\times 10^{-4}$ & $2.98\times 10^{-2}$ \\
 & Control error (MSE) & $6.11\times 10^{-2}$ & $3.71\times 10^{-2}$ & $3.54\times 10^{-2}$ & $2.41\times 10^{-1}$ & $8.3\times 10^{-2}$ & $2.78\times 10^{-2}$ & $6.43\times 10^{-2}$ \\
\bottomrule
\end{tabular}
\end{table}

\subsection{Nonlinear Transport Control: Viscous Burgers' Equation}

\begin{figure}
    \centering

    \begin{subfigure}{\linewidth}
        \centering
        \includegraphics[width=0.90\linewidth]{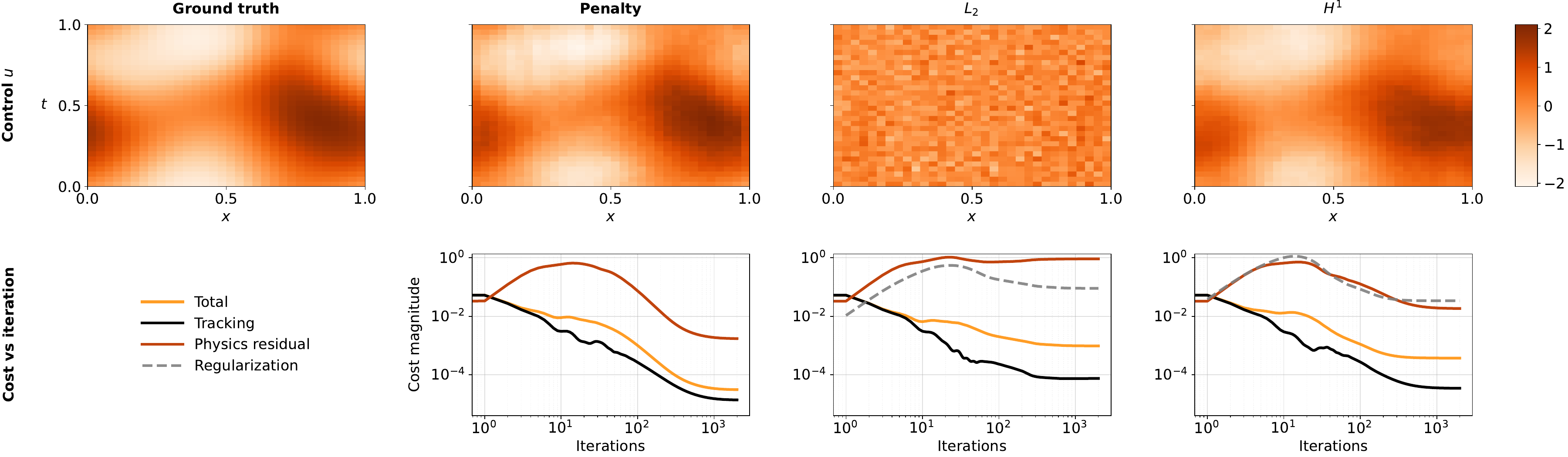}
        \caption{Objective $J_{Burg,1}$.}

        \label{fig:burgers_1_regularization}
    \end{subfigure}

    \vspace{1pt} 

    \begin{subfigure}{\linewidth}
        \centering
        \includegraphics[width=0.90\linewidth]{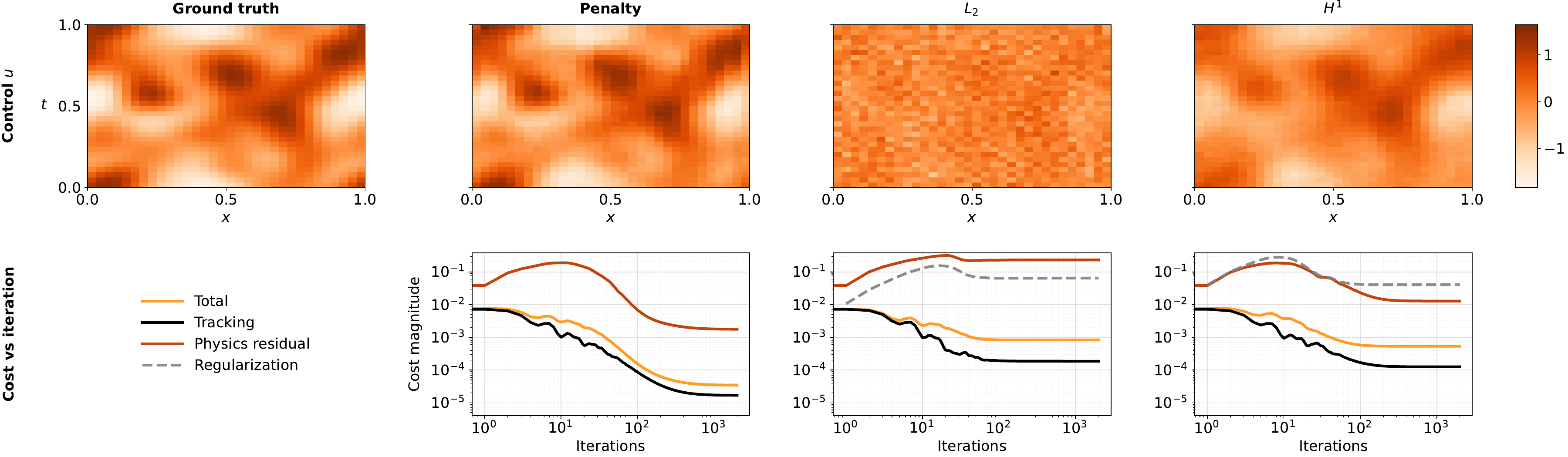}
        \caption{Objective $J_{Burg,2}$.}
        \label{fig:burgers_2_regularization}
    \end{subfigure}

    \caption{Nonlinear transport control problem \eqref{eq:burgers} with costs \eqref{cost:burgers}: comparison of three NO-approaches against the ground truth. Columns (left to right) show: ground truth (reference), \emph{penalty-only} method ($\lambda=0$), \emph{$L_2$-regularization only} ($\mu=0$), and \emph{$H^1$-regularization only} ($\mu=0$). Top row shows final control $u$ and bottom row optimization cost versus iteration.}
    \label{fig:burgers_regularization}
\end{figure}

\begin{figure}
    \centering

    \begin{subfigure}{\linewidth}
        \centering
        \includegraphics[width=0.85\linewidth]{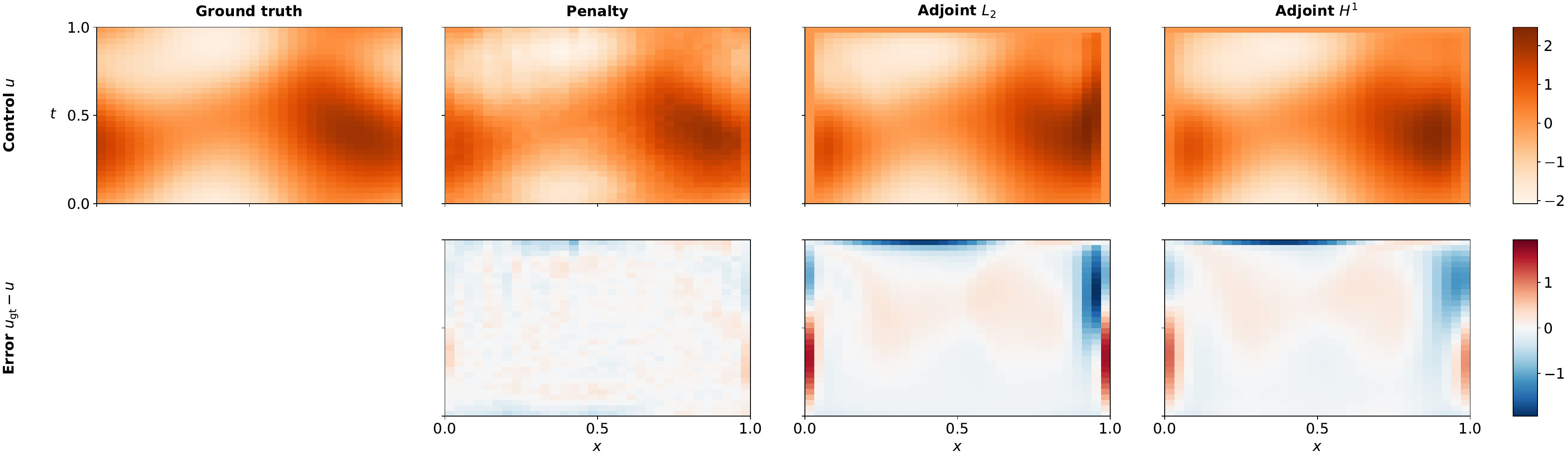}
        \caption{Objective $J_{Burg,1}$.}
        \label{fig:burgers_1_regularization_reference}
    \end{subfigure}

    \vspace{1pt} 

    \begin{subfigure}{\linewidth}
        \centering
        \includegraphics[width=0.85\linewidth]{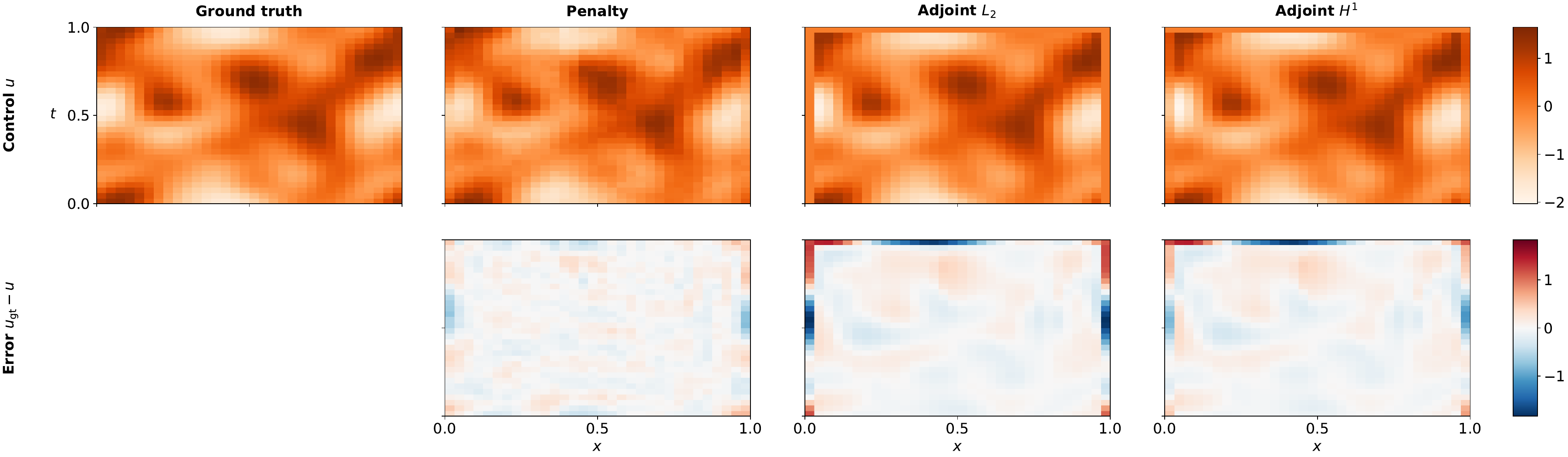}
        \caption{Objective $J_{Burg,2}$.}
        \label{fig:burgers_2_regularization_reference}
    \end{subfigure}

\caption{Nonlinear transport control problem \eqref{eq:burgers} with costs \eqref{cost:burgers}: comparison of three approaches against the ground truth. Columns (left to right) show: ground truth (reference), \emph{penalty-only} method ($\lambda=0$), adjoint method with \emph{$L_2$-regularization}, and adjoint method with \emph{$H^1$-regularization}. Top row shows final control $u$ and bottom row control error $u-u_{ref}$.}
    \label{fig:burgers_regularization_reference}
\end{figure}

Next, we consider the tracking problem \eqref{cost:burgers} and with the nonlinear Burgers' equation \eqref{eq:burgers} acting as the constraint. Similarly to the scalar elliptic control problem, we study the effects of the penalty and the regularization for the problems $J_{Burg,1}$ and $J_{Burg,2}$. Again, we demonstrate the existence and effects of the cheating directions in Figures \ref{fig:burgers_1_regularization} and \ref{fig:burgers_2_regularization}, where $L_2$ regularization fails to find a feasible solution. Similarly, $H^1$ finds a more feasible solution and smooth solution, but is outperformed by the penalty method in terms of feasibility. 

In contrast to the scalar elliptic control problem, this problem is time-dependent. Hence, the adjoint method requires a time-stepping scheme to calculate the forward solution as well as the gradient through the adjoint equation. These time-stepping schemes are more error-prone as errors can accumulate over time. We observe this effect with higher interior error as can be seen in Figures \ref{fig:burgers_1_regularization_reference} and \ref{fig:burgers_2_regularization_reference}. In this regard, the NO shows superior performance. Similarly to the scalar elliptic control problem, the effect of Dirichlet boundary conditions can be seen in the error of the adjoint method, i.e, the control points at the boundary do not influence the solution.

\begin{table}[t]
\centering
\scriptsize
\caption{Summary of errors for nonlinear transport control problem for different penalties, $L_2$, and $H^1$ regularization.}
\label{tab:summary_errors_burgers}
\renewcommand{\arraystretch}{1.15}
\setlength{\tabcolsep}{3pt}

\begin{tabular}{c l ccc c ccc}
\toprule
& & \multicolumn{3}{c}{Penalty method} & \multicolumn{1}{c}{$L_2$ regularization} & \multicolumn{3}{c}{$H^1$ regularization} \\
\cmidrule(lr){3-5}\cmidrule(lr){6-6}\cmidrule(lr){7-9}
Problem & Metric & $\mu=10^{-1}$ & $\mu=10^{-2}$ & $\mu=10^{-3}$ & $\lambda=10^{-2}$ & $\lambda=10^{-1}$ & $\lambda=10^{-2}$ & $\lambda=10^{-3}$ \\
\midrule
\multirow{3}{*}{$J_{Burg,1}$} & Tracking error $\mathcal{J}_{\text{track}}$ & $3.14\times 10^{-5}$ & $1.38\times 10^{-5}$ & $1.42\times 10^{-5}$ & $7.51\times 10^{-5}$ & $5.17\times 10^{-4}$ & $3.48\times 10^{-5}$ & $1.43\times 10^{-5}$ \\
 & Physics residual $\mathcal{R}_\theta$ & $3.46\times 10^{-4}$ & $1.73\times 10^{-3}$ & $1.14\times 10^{-1}$ & $9.15\times 10^{-1}$ & $1.23\times 10^{-2}$ & $1.84\times 10^{-2}$ & $1.87\times 10^{-1}$ \\
 & Control error (MSE) & $8.91\times 10^{-3}$ & $1.49\times 10^{-2}$ & $1.64\times 10^{-1}$ & $9.68\times 10^{-1}$ & $1.01\times 10^{-1}$ & $4.05\times 10^{-2}$ & $2.27\times 10^{-1}$ \\
\midrule
\multirow{3}{*}{$J_{Burg,2}$} & Tracking error $\mathcal{J}_{\text{track}}$ & $5.35\times 10^{-5}$ & $1.7\times 10^{-5}$ & $1.19\times 10^{-5}$ & $1.87\times 10^{-4}$ & $1.13\times 10^{-3}$ & $1.26\times 10^{-4}$ & $1.01\times 10^{-5}$ \\
 & Physics residual $\mathcal{R}_\theta$ & $5.77\times 10^{-4}$ & $1.77\times 10^{-3}$ & $1.09\times 10^{-2}$ & $2.33\times 10^{-1}$ & $3.91\times 10^{-3}$ & $1.29\times 10^{-2}$ & $2.09\times 10^{-1}$ \\
 & Control error (MSE) & $2.66\times 10^{-2}$ & $1.99\times 10^{-2}$ & $5.05\times 10^{-2}$ & $4.26\times 10^{-1}$ & $2.08\times 10^{-1}$ & $7.68\times 10^{-2}$ & $2.58\times 10^{-1}$ \\
\bottomrule
\end{tabular}
\end{table}


\subsection{Flow Control: The Stokes Equation}
In our third experiment, we consider a flow control (inverse) problem \eqref{cost:stokes}. Given velocity measurements, we infer the unknown forcing field—treated as the control disturbing the flow, which is governed by the Stokes equations \eqref{eq:stokes}. In contrast to the scalar elliptic control and the nonlinear transport control problems, the problem is ill-posed in the sense that it is underdetermined, as we are not tracking the pressure field, i.e., the velocity field is not unique, and multiple controls can generate the same velocity field when the pressure is not fixed. Therefore, we do not do a posterior error analysis on the control but rather study the feasibility of the optimized control through a reference solver. In other words, we first optimize and get a candidate solution $\mathbf{u}_{S,\mu}$, which we then use to solve the velocities $\mathbf{v}(\mathbf{u}_{S,\mu})$ with our reference forward PDE solver and compare the error in the velocities. 

\begin{figure}  
    \centering
    \includegraphics[width=0.99\linewidth]{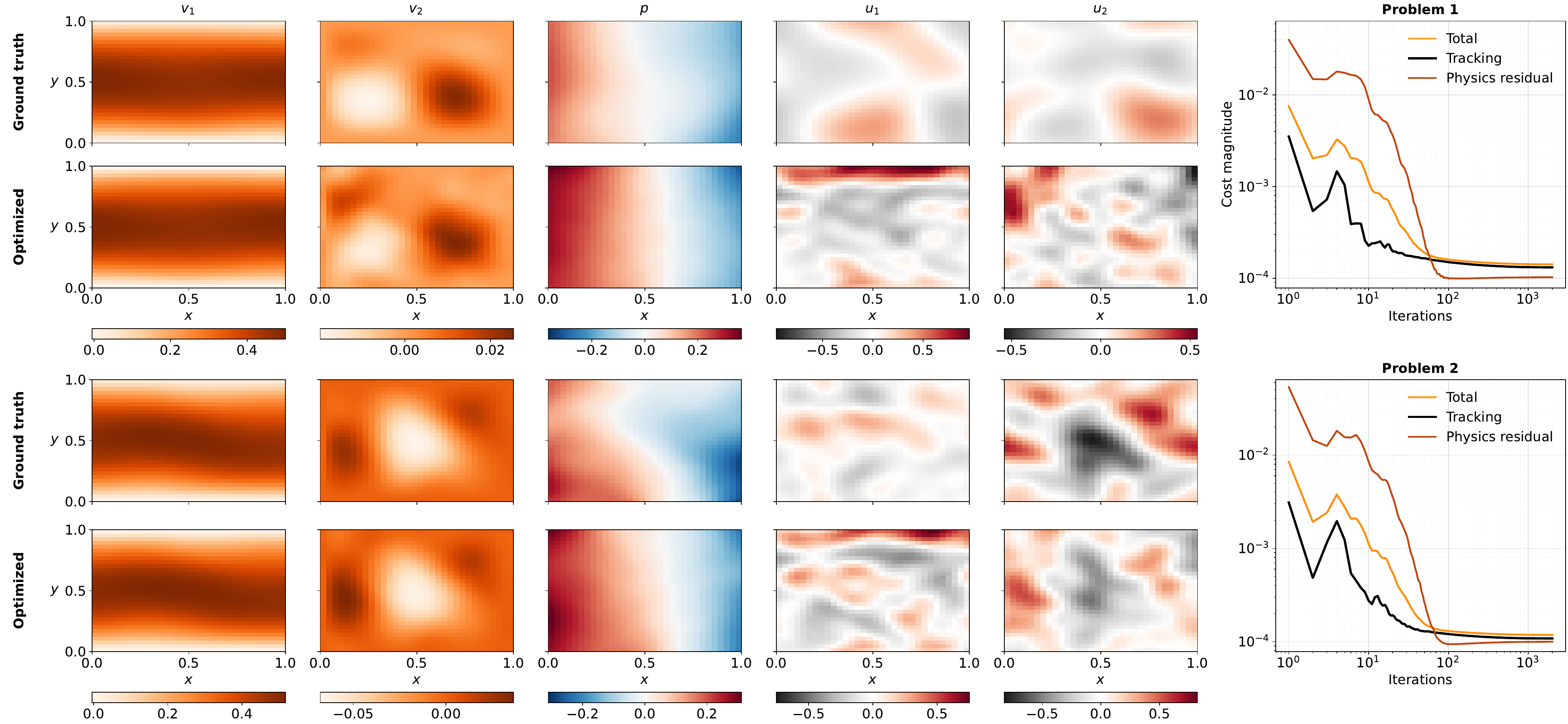}  
\caption{Flow control problem \eqref{eq:stokes} and \eqref{cost:stokes}: comparison of the optimized solution against the ground truth for Problem 1 (top half) and Problem 2 (bottom half). Columns (left to right) show: velocity components $v_1$ and $v_2$, pressure $p$, and control fields $u_1$ and $u_2$. For each problem, the top row displays the ground truth result and the bottom row the optimized solution. The rightmost column plots the optimization cost components (Total, Tracking, and Physics residual) versus iterations.}
    \label{fig:stokes_convergence}
\end{figure}

In Figure \ref{fig:stokes_convergence}, we can see that the controls differ from the reference control as well as the pressure field. However, the problem is ill-posed in the sense that the solution is not unique, as we are not tracking the pressure. From the convergence graph, we observe that the differential equation is not violated as the residual is minimized to $10^{-4}$. To verify that the optimized control truly generates the velocity field $\mathbf{v}_d$, we compare the generated state when given the optimized control with a reference finite difference solver. The details of the reference solver can be found in \ref{appendix_b}. From Figure \ref{fig:stokes_FD_both}, we see that the optimized control generates the target velocity field $\mathbf{v}_d$.

\begin{figure}
    \centering

    \begin{subfigure}{\linewidth}
        \centering        
        \includegraphics[width=0.99\linewidth]{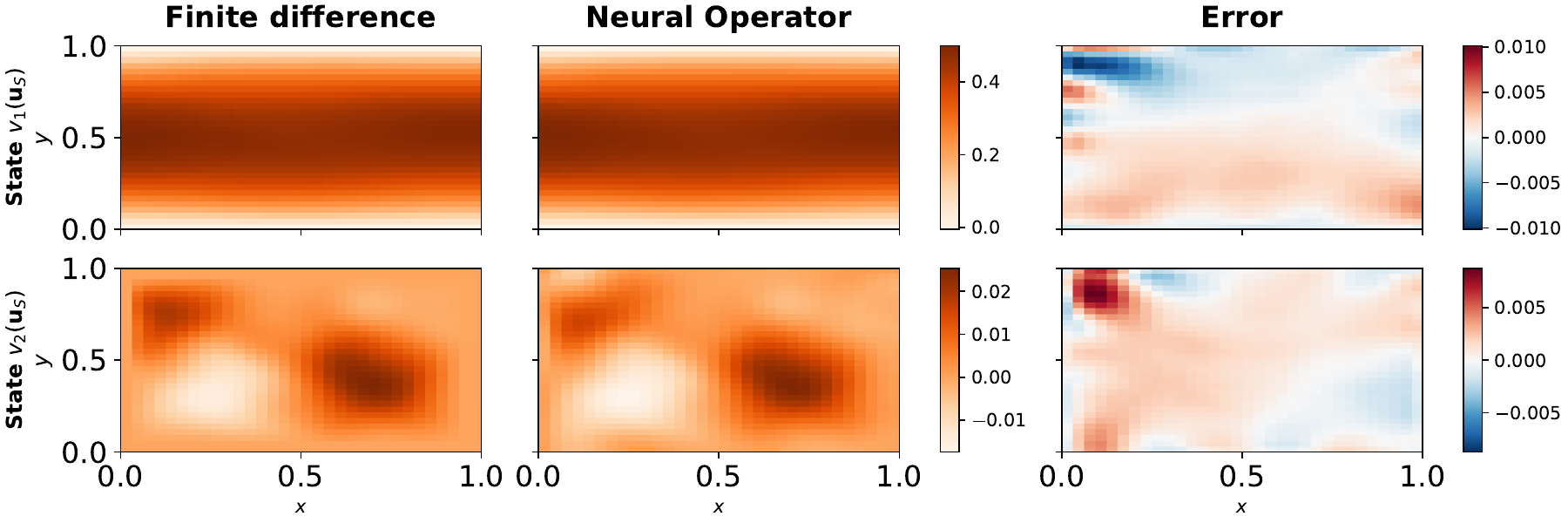}
        \caption{$J_{Stokes,1}$}
        \label{fig:stokes_FD_1}
    \end{subfigure}

    \vspace{1pt} 

    \begin{subfigure}{\linewidth}
        \centering     
        \includegraphics[width=0.99\linewidth]{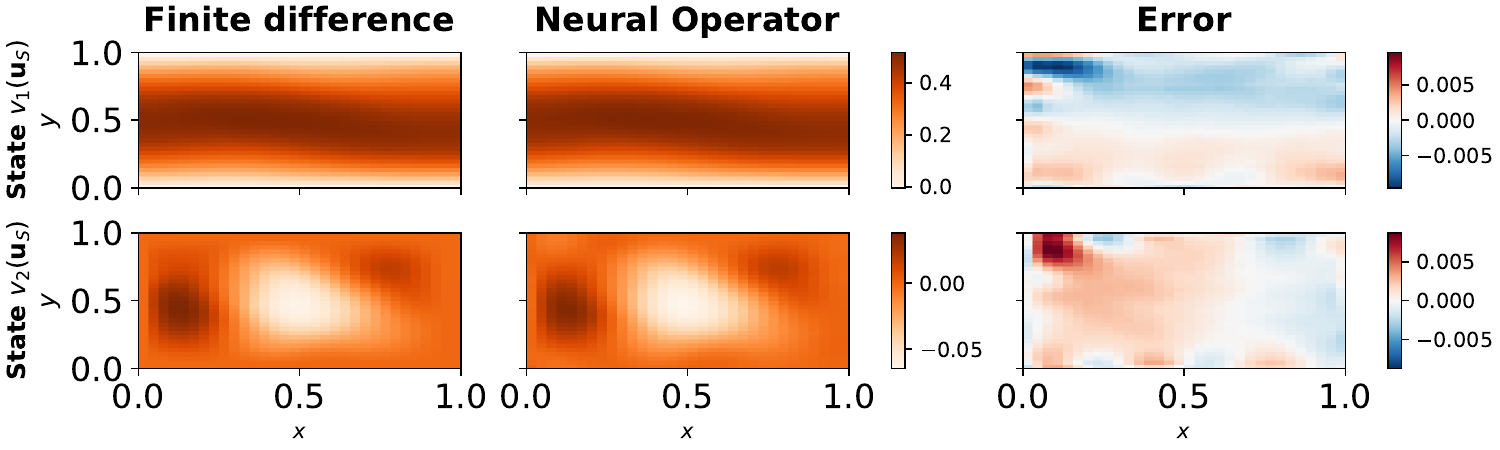}
        \caption{$J_{Stokes,2}$}
        \label{fig:stokes_FD_2}
    \end{subfigure}

    \caption{Velocity states $v_1$ (top) and $v_2$ (bottom) generated by the optimized control $\mathbf{u}_{S,\mu}$ verified by a finite difference solver (left column) and the DeepONet (middle column) and their difference (right column).}
    \label{fig:stokes_FD_both}
\end{figure}

\subsection{Sensitivity Analysis}
\begin{figure}[!t]
    \centering

    \begin{subfigure}{\linewidth}
        \centering
        \includegraphics[width=0.8\linewidth]{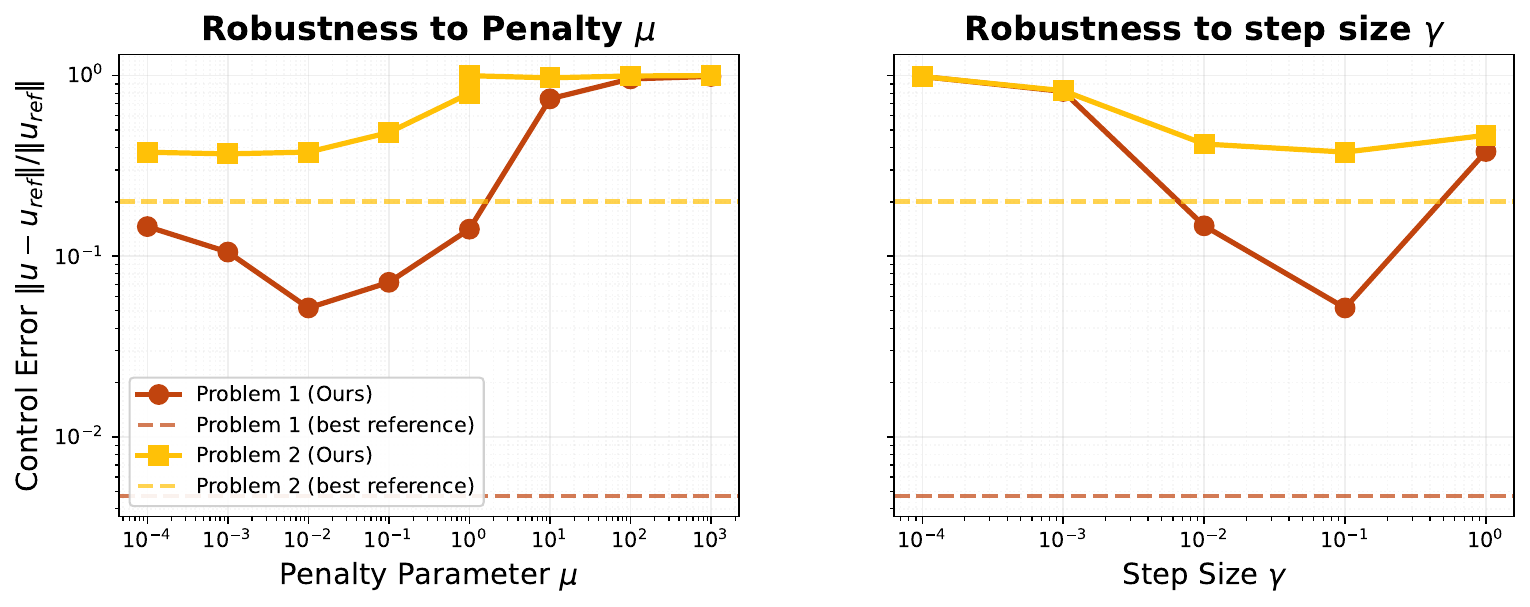}
        \caption{Poisson's problem.}
        \label{fig:poisson_sensitivity}
    \end{subfigure}

    \vspace{1pt} 

    \begin{subfigure}{\linewidth}
        \centering
        \includegraphics[width=0.8\linewidth]{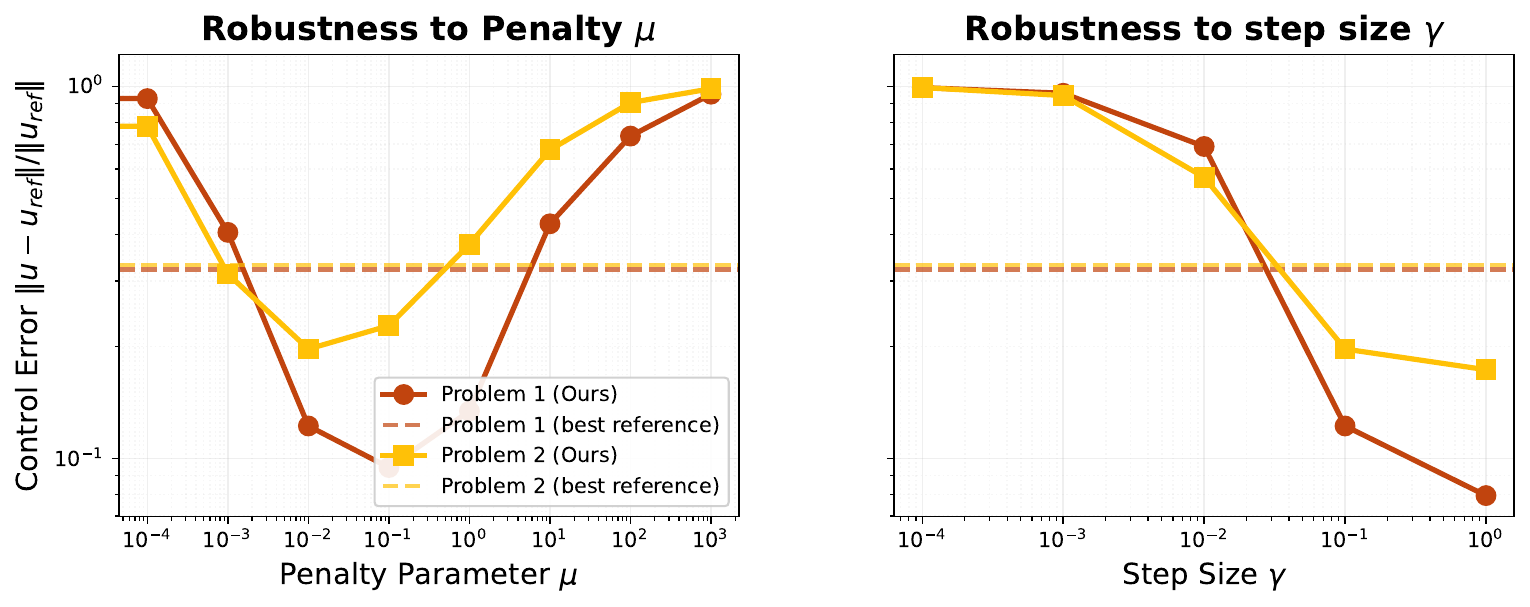}
        \caption{Nonlinear transport control problem.}
        \label{fig:burgers_sensitivity}
    \end{subfigure}

    \caption{Sensitivity anaylsis for the scalar elliptic control and the nonlinear transport control problem.  Left is the  relative error vs. penalty parameter $\mu$ and right is the relative error vs. gradient step size $\gamma$.}
    \label{fig:sensitivity_both}
\end{figure}

For the scalar elliptic control and the nonlinear transport control problems, we study the effects of the penalty factor $\mu$ and the initial step size $\gamma$. We perform a sweep of the parameters in the ranges $\mu\in[10^{-4},10^3]$ and $\gamma=[10^{-4},1]$ by increasing the parameter tenfold for each run. We compare the posterior relative error $\|\mathbf{u}_S-\mathbf{u}_{ref}\|/\|\mathbf{u}_{ref}\|$. For comparison, we add the best result (smallest error) for the adjoint reference method with either $L_2$ or $H^1$ regularizer as a dashed line in the figure. From Figures \ref{fig:poisson_sensitivity} and \ref{fig:burgers_sensitivity}, we can conclude that the method is only mildly sensitive to the penalty parameter, but it should be of the right magnitude for best performance. 

For the initial step size $\gamma$ of the gradient update, we observe that a larger, but still reasonable, step size is better. Thus, an initial step size in the range of $[0.1,0.5]$ provides good general performance. The low sensitivity to the step size is due to the decay of the step size with the AdamW optimizer. This, however, opens up the possibility of reducing the number of iterations to improve computational speed, which was not the main focus in this study. Instead, we used a fixed number of iterations for all cases for simplicity and transparent comparison.

We did not study the sensitivity of the optimization to the training quality of the NO. In practice, physics-informed training for complex PDEs may stall at higher residual levels than those achieved here. However, we note that the residual penalty term in the objective~\eqref{eq:pde_tracking_no_penalized_reduced} provides an implicit safeguard: a poorly trained NO will produce large PDE residuals, which inflate the penalized objective and prevent the optimizer from accepting controls whose corresponding states are inaccurately approximated. Thus, the residual penalty not only enforces feasibility but also acts as a self-diagnostic for surrogate quality during optimization. A systematic study of how optimization accuracy degrades as a function of training residual is left to future work.

\subsection{Solution Metrics}
We summarize solution times in Table \ref{tab:runtime_metrics_pois_burg} and accuracy in terms of MSE in Table \ref{tab:accuracy_metrics_pois_burg} for the scalar elliptic control and the nonlinear transport control problems. For the scalar elliptic control problem, the neural operator also performs slightly worse than the adjoint methods but remains competitive in both accuracy and solution time. The reason for the adjoint method is that Poisson's equation is a linear PDE, which can thus be solved accurately and quickly. We also compare the accuracy at interior points to disregard boundary effects.

For the nonlinear transport control problem, our method achieves substantially faster solution times, with about four times faster and with better accuracy. This advantage arises because the PDE is time-dependent; thus, each optimization step with the adjoint method requires marching forward in time across all spatial coordinates, which is less efficient than solving a linear system.

Our method requires an upfront training cost of approximately 90 minutes for the scalar elliptic control (Poisson) and nonlinear transport (Burgers) problems. However, this cost is incurred once: the trained NO can be reused for any number of downstream control problems on the same PDE class. For the nonlinear transport control problem, each surrogate optimization solve requires approximately 16 seconds compared to 73 seconds for the adjoint method, yielding a marginal saving of roughly 57 seconds per problem. The upfront training cost is therefore amortized after approximately 95 control problems, which comprises a realistic scenario in engineering applications involving repeated design iterations, parameter studies, or real-time re-optimization with changing targets. This stands in contrast to PINN-based and multi-network approaches, where the training cost is incurred for every new cost function or tracking target, scaling linearly with the number of problems solved.

\begin{table*}[t]
\centering
\caption{Runtime metrics for the scalar elliptic control and the nonlinear transport control problems. Best (lower) per row and metric group in \textbf{bold}.}
\label{tab:runtime_metrics_pois_burg}
\scriptsize
\setlength{\tabcolsep}{3pt}
\begin{tabular}{lcccccc}
\toprule
\textbf{Problem} & \multicolumn{3}{c}{\textbf{Solution time [s]}} & \multicolumn{3}{c}{\textbf{Time/iteration [s]}} \\
\cmidrule(lr){2-4}\cmidrule(lr){5-7}
 & \textbf{Ours} & \textbf{Adj.\ $L^2$} & \textbf{Adj.\ $H^1$} & \textbf{Ours} & \textbf{Adj.\ $L^2$} & \textbf{Adj.\ $H^1$} \\
\midrule
$J_{\mathrm{Pois},1}(u)$ & $\boldsymbol{8.67}$ & $9.58$ & $9.94$ & $\boldsymbol{4.33\times 10^{-3}}$ & $4.79\times 10^{-3}$ & $4.97\times 10^{-3}$ \\
$J_{\mathrm{Pois},2}(u)$ & $\boldsymbol{8.64}$ & $13.6$ & $13.4$ & $\boldsymbol{4.32\times 10^{-3}}$ & $6.8\times 10^{-3}$ & $6.69\times 10^{-3}$ \\
$J_{\mathrm{Burg},1}(u)$ & $\boldsymbol{16.1}$ & $74.8$ & $72$ & $\boldsymbol{8.07\times 10^{-3}}$ & $0.0374$ & $0.036$ \\
$J_{\mathrm{Burg},2}(u)$ & $\boldsymbol{16.3}$ & $75.7$ & $72.2$ & $\boldsymbol{8.13\times 10^{-3}}$ & $0.0379$ & $0.0361$ \\
\bottomrule
\end{tabular}
\end{table*}

\begin{table*}[t]
\centering
\caption{Accuracy metrics for the scalar elliptic control and the nonlinear transport control problems. Best (lower) per row and metric group in \textbf{bold}.}
\label{tab:accuracy_metrics_pois_burg}
\scriptsize
\setlength{\tabcolsep}{2pt}
\begin{tabular}{lccccccccc}
\toprule
\textbf{Problem} & \multicolumn{3}{c}{\textbf{Control error (all)}} & \multicolumn{3}{c}{\textbf{Control error (interior, excl.\ boundary)}} & \multicolumn{3}{c}{\textbf{Tracking error}} \\
\cmidrule(lr){2-4}\cmidrule(lr){5-7}\cmidrule(lr){8-10}
 & \textbf{Ours} & \textbf{Adj.\ $L^2$} & \textbf{Adj.\ $H^1$} & \textbf{Ours} & \textbf{Adj.\ $L^2$} & \textbf{Adj.\ $H^1$} & \textbf{Ours} & \textbf{Adj.\ $L^2$} & \textbf{Adj.\ $H^1$} \\
\midrule
$J_{\mathrm{Pois},1}(u)$ & $6.29\times 10^{-4}$ & $\boldsymbol{5.26\times 10^{-6}}$ & $3.62\times 10^{-4}$ & $6.24\times 10^{-4}$ & $\boldsymbol{5.99\times 10^{-6}}$ & $8.66\times 10^{-5}$ & $5.36\times 10^{-9}$ & $9.63\times 10^{-9}$ & $\boldsymbol{1.05\times 10^{-10}}$ \\
$J_{\mathrm{Pois},2}(u)$ & $0.0371$ & $0.0218$ & $\boldsymbol{0.0105}$ & $0.0221$ & $\boldsymbol{2.78\times 10^{-4}}$ & $1.69\times 10^{-3}$ & $7.58\times 10^{-8}$ & $\boldsymbol{8.19\times 10^{-11}}$ & $9.39\times 10^{-10}$ \\
$J_{\mathrm{Burg},1}(u)$ & $\boldsymbol{0.0149}$ & $0.172$ & $0.104$ & $\boldsymbol{7.57\times 10^{-3}}$ & $0.0909$ & $0.0476$ & $\boldsymbol{1.39\times 10^{-5}}$ & $1.01\times 10^{-4}$ & $1.62\times 10^{-4}$ \\
$J_{\mathrm{Burg},2}(u)$ & $\boldsymbol{0.0199}$ & $0.0984$ & $0.0565$ & $\boldsymbol{9.44\times 10^{-3}}$ & $0.0126$ & $0.0132$ & $1.7\times 10^{-5}$ & $7.77\times 10^{-7}$ & $\boldsymbol{3.87\times 10^{-9}}$ \\
\bottomrule
\end{tabular}
\end{table*}

\section{CONCLUSIONS}\label{section_6}

We show that a physics-informed DeepONet model can be directly applied to solve control problems of tracking type or control recovery together with an unconstrained optimizer, without requiring any information about the cost functions during training. This significantly enhances the practical applicability of NOs and reduces the architectural complexity typically associated with solving control problems using neural networks. Our results further show that neural operators can be used beyond their original role as differential equation solvers. Unlike conventional neural network approaches, we did not need to employ dedicated control networks or auxiliary components to explicitly construct the solution. Instead, the trained DeepONet was integrated into an optimization routine, where the optimal control problem was discretized and a differential-equation residual was added to the cost as a penalty, ensuring that the control remained within the solution space during iterative optimization. We discussed and showed how the residual penalty effectively acts as a low-pass filter for PDEs with a damping or dissipative term, eliminating the need for a regularizer to enforce well-posedness.

For the scalar elliptic control problem, our method showed inferior control accuracy compared to the reference adjoint method, which we argue is due to the linear property of the problem. Thus, solving the scalar elliptic control problem with the adjoint method reduces to solving linear transformations while updating the gradient, which can be done effectively and to high accuracy. However, for the nonlinear transport control, which has a time-dependent PDE (Burgers’ equation), our approach yielded iteration times up to four times faster than the reference adjoint method, highlighting its potential for more complex time-dependent PDE-control problems. A clear advantage of our approach is that, for PDE-control, we never need to implement a PDE solver or solution scheme within the optimization loop itself; instead, the physics-informed training of the neural operator embeds the dynamics directly.

While our study demonstrated the potential of physics-informed DeepONets for optimal control, we focused on a specific set of assumptions and an open-loop setting. Our investigation was limited to the DeepONet architecture, albeit with some variations in its formulation. However, we see no limitations regarding the use of alternative NO architectures with likely similar or better performance, as long as they can be trained as physics-informed. This, however, limits the method's usability, since training an NO in a purely physics-informed way is not trivial for complex differential equations. In addition, if the input is not normalized or within a suitable range, physics-informed training becomes difficult, further restricting general applicability. We also restricted the class of admissible controls to smooth functions, as NOs cannot effectively handle discontinuous functions by default. This further simplifies the analysis and optimization, but may not capture discontinuous or bang-bang control strategies. Finally, we assumed bounded solution spaces, reachable states and recoverable controls, leaving the treatment of more complex or unbounded systems, as well as more complicated PDEs and unreachable states, for future work. Further, we did not consider terminal costs, but our preliminary experiments, as shown in \ref{appendix_d}, suggest that the method works for terminal and quadratic costs as well, when an additional smoothing regularization is used together with the residual penalty. A detailed explanation for why the approach remains effective under such additional regularization is beyond the scope of this work and is left for future investigation. 

Overall, we view our encouraging results as a first step toward implementing physics-informed NO as such, e.g, without additional architectural modifications, as surrogate models for control and inverse problems. Lastly, as for future studies, we see providing a more theoretical framework to study when a physics-informed neural operator, such as the DeepONet can act as an efficient surrogate for PDE-control problems, and for what types of problems, as likely fruitful avenues to pursue.

\section*{ACKNOWLEDGMENT}

This work was supported by the Finnish Ministry of Education and Culture’s Pilot for Doctoral Programmes (Pilot project Mathematics of Sensing, Imaging and Modelling).

\newpage
\appendix

\section{Architectures and Training}\label{appendix_a}

We used a DeepONet model as a base model for our NO, with a modified architecture instead of the fully connected architecture. The modified network has a skip connection for each layer to improve the training of the NO. The architecture for the modified network is presented by Wang et al.\cite{wang2021understanding}, who refer to it as as an "improved fully connected network". The models are trainable with FCN and the method works for the FCN, but we found the models to be more easy to train with modified network of Wang et al. \cite{wang2021understanding}.

The dimension sizes of the layers are shown in Table \ref{tab:appendix_architecture}. For all models we used a hyperbolic tangent activation function. The final output dimensions of the trunk network is 1024, for it to match with the dot product with the branch network. We used no bias parameters at the final dot product between the trunk and the branch. 

\begin{table}[h!]
\centering
\caption{DeepONet model configurations.}
\label{tab:appendix_architecture}
\begin{tabular}{lccc}
\toprule
\textbf{Model} & $m$  & \makecell{Branch\\(W $\times$ D)} & \makecell{Trunk\\(W $\times$ D)} \\
\midrule
Poisson's eq.             & $[32,32]$     & $1024 \times 3$ & $200 \times 5$ \\
Burgers' eq. & $[32,32]$     & $1024 \times 3$ & $400 \times 6$ \\
Stoke's eq.             & $[32,32]$     & $1024 \times 3$ & $300 \times 4$ \\
\bottomrule
\end{tabular}
\end{table}

The NOs are trained on different input function sets where the parameters for the functions are sampled uniformly from a given interval, on their respective domain. We used polynomial functions and Gaussian Random Fields (GRF), defined as
\begin{equation}\label{function_sets}
\begin{aligned}
    \text{Polynomial}:&\quad a_n t^n + a_{n-1} t^{n-1} + \dots + a_1 t + a_0\\
    \text{GRF}:&\quad \mathcal{GP}(0, \sigma^2 \exp ( -\frac{|t - t'|^2}{2l^2})).
\end{aligned}
\end{equation}

The training set is constructed of an equal number of functions drawn from a chosen subset of the types listed in \eqref{function_sets}. In addition, the parameters for each function are sampled from a uniform distribution within the specific ranges. The subsets for each model and the ranges for the parameter sampling are given in Table \ref{tab:function_sets_parameter}. The parameter ranges were chosen to ensure that the solution of the differential equation does not exhibit excessively rapid growth or decay. Further, the input functions $u$, where projected such that their max or min value was within a reasonable domain for the problem. We use 300,000 functions for training the models. We used a batch size of 100 for the functions (branch input) as well as 100 sampled time, spatial coordinates or time-spatial coordinates (trunk input). This corresponds to $300,000/100=3000$ epochs in conventional machine learning terms. We trained the models purely as physics-informed, i.e., unsupervised, only based on the residual of the differential equations with equal weighting. We used Optax \cite{deepmind2020jax} ADAM optimizer with the default parameters of $\beta=(0.9, 0.999)$ and a learning rate scheduler that decreases the learning rate by a decade after $100,000$ steps. 

\begin{figure}[t]
    \centering
    \includegraphics[width=0.99\linewidth]{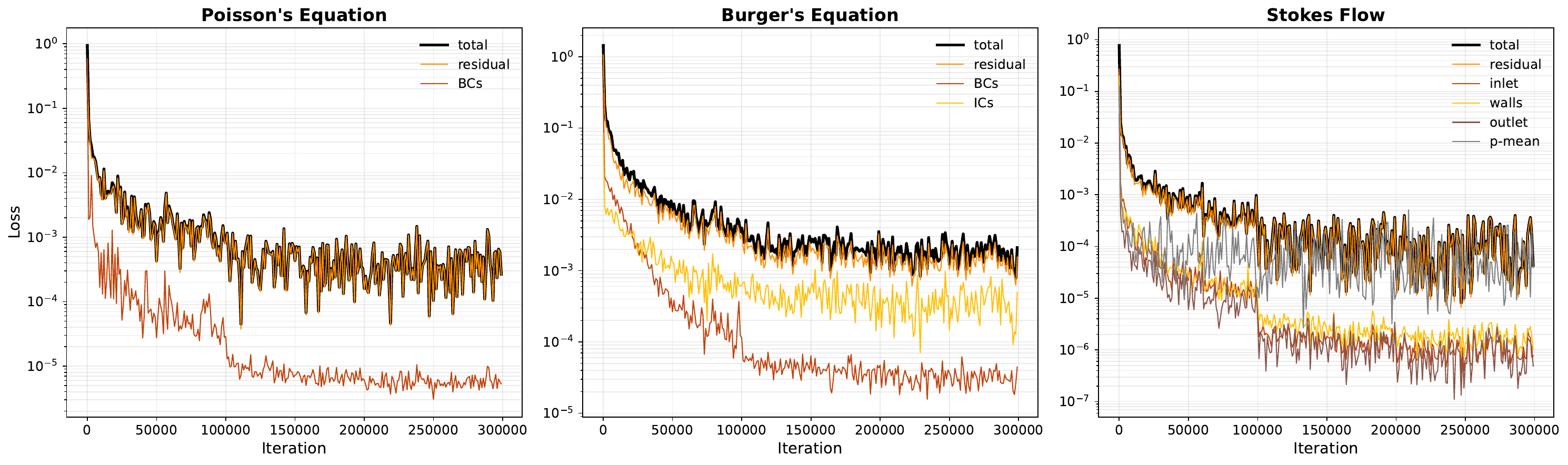}  
    \caption{Training loss for the DeepONet models.}
    \label{fig:training_metrics}
\end{figure}

\begin{table}
    \caption{Function sets and parameter values for different problems.}
    \label{tab:function_sets_parameter}
    \centering
    \renewcommand{\arraystretch}{1.0}
    \begin{tabular}{l|c c c c}
        \toprule
        \multirow{2}{*}{\textbf{Model}} & \multicolumn{1}{c}{\textbf{GRF}} & \multicolumn{2}{c}{\textbf{Polynomial}} & \multicolumn{1}{c}{\textbf{Function min/max}}\\
        \cmidrule(lr){2-2} \cmidrule(lr){3-4} \cmidrule(lr){5-5}
        &  $ l $&  $n$ & $a_n$ & $[u_{min},u_{max}]$\\
        \midrule
        Poisson's eq. & $ [0.2, 2.0] $ & $\{0,\dots,8\}$ & [-2.0,2.0] & $[-2,2]$ \\
        Forced Burgers' eq. & $ [0.2, 3.0] $ & $\{0,\dots,8\}$ & [-2.0,2.0] & $[-2,2]$ \\
        Stokes & $ [0.2, 2.0] $ & $\{0,\dots,8\}$ & [-2.0,2.0] & $[-1.2,1.2]$ \\
        \bottomrule
    \end{tabular}
    \label{tab:function_sets}
\end{table}

We used JAX \cite{jax2018github} for the implementation. The training was performed on a standard desktop computer running an Intel Core i5-12600K processor and NVIDIA RTX A2000 GPU. The training times were approximately 90 minutes for Poisson and Burgers,  and 150 minutes for Stokes flow.

The training loss for the total (sum of all), differential equation residual, initial and boundary losses are shown in Figure \ref{fig:training_metrics}. We do not monitor any test or validation set as we are not directly interested in how well the network generalizes on unseen data, but rather in whether it can solve the optimal control problem. In addition, as our training function set is quite large, creating a test set that does not include any function of the training set is cumbersome and would most likely require calculating correlations between each function in the test and train sets. Hence, we deemed it to be not relevant enough to justify the effort, as the final tests are done in the optimal control framework.   

\clearpage

\section{Reference Methods}\label{appendix_b}
\subsection{Scalar Elliptic Control: The Poisson Equation}

For the adjoint method, the state and adjoint equations for the Poisson tracking problems \eqref{cost:poisson} are obtained from the first-order optimality conditions. We recall the Poisson state equation
\begin{equation}\label{eq_app:poisson_state}
\begin{aligned}
\Delta y(\mathbf{x}) &= -u(\mathbf{x}) \qquad \mathbf{x\in} \Omega,\\
y(\mathbf{x}) &= 0 \qquad \qquad \mathbf{x} \in \partial\Omega,
\end{aligned}
\end{equation}
and its adjoint equation for the tracking term,
\begin{equation}\label{eq_app:poisson_adjoint}
\begin{aligned}
\Delta p(\mathbf{x}) &= y(\mathbf{x}) - y_d(\mathbf{x}) \qquad \mathbf{x\in} \Omega,\\
p(\mathbf{x}) &= 0 \qquad \qquad \qquad \mathbf{x} \in \partial\Omega,
\end{aligned}
\end{equation}
The pair \eqref{eq_app:poisson_state}--\eqref{eq_app:poisson_adjoint} forms a linear system that can be solved sequentially: given a control $u_k$, one first solves for the state $y_k$, and then solves for the adjoint $p_k$.

To update the control, we compute the gradient of the Lagrangian with respect to $u$.
For an $L_2$ control regularization term $\lambda\|u\|_{2}^2$, the gradient is
\begin{equation}\label{eq_app:poisson_adjoint_L2_grad}
\nabla_u \mathcal{L}(u,y,p)=2\lambda u(\mathbf{x}) - p(\mathbf{x}).
\end{equation}
For an $H^1$-seminorm regularization term $\lambda\|\nabla u\|_{2}^2$, the gradient takes the form
\begin{equation}\label{eq_app:poisson_adjoint_H1_grad}
\nabla_u \mathcal{L}(u,y,p)= -2\lambda \Delta u(\mathbf{x}) - p(\mathbf{x}),
\end{equation}
with the corresponding natural boundary condition $\partial_n u=0$ on $\partial\Omega$ if $u$ is not prescribed on the boundary.

We perform a gradient descent step,
\begin{equation}\label{eq_app:poisson_gradient_descent}
u_{k+1}(\mathbf{x}) = u_k(\mathbf{x}) - \gamma\,\nabla_u \mathcal{L}(u_k,y_k,p_k),
\end{equation}
where $k$ denotes the iteration index and $\gamma>0$ is the step size. In our JAX implementation, each iteration consists of:
\begin{enumerate}
    \item Solve the state $y_k$ from \eqref{eq_app:poisson_state} for the current control $u_k$.
    \item Solve the adjoint $p_k$ from \eqref{eq_app:poisson_adjoint} using the state $y_k$.
    \item Evaluate the gradient $\nabla_u \mathcal{L}(u_k,y_k,p_k)$ via \eqref{eq_app:poisson_adjoint_L2_grad} or \eqref{eq_app:poisson_adjoint_H1_grad}.
    \item Update the control using \eqref{eq_app:poisson_gradient_descent}.
\end{enumerate}

To solve the state equation \eqref{eq_app:poisson_state} and \eqref{eq_app:poisson_adjoint}, we use a finite difference scheme and approximate the Laplacian operator with a second-order approximation (5-point-stencil). Let $x_i = ih$, $y_j = jh$, $y_{i,j} \approx y(x_i,y_j)$ and $h$ be the grid spacing. Then, for $1\le i,j \le N-2$,
\begin{equation}\label{eq:laplacian_5pt}
(\Delta_h y)_{i,j}
:= \frac{y_{i+1,j} + y_{i-1,j} + y_{i,j+1} + y_{i,j-1} - 4y_{i,j}}{h^2},
\end{equation}
which satisfies
\begin{equation}\label{eq:laplacian_5pt_error}
(\Delta y)(x_i,y_j) = (\Delta_h y)_{i,j} + \mathcal{O}(h^2).
\end{equation}
We use the same grid of sensor locations $S$ as for the neural operator, i.e., a uniform $(32,32)$ grid. To obtain a symmetric positive definite linear system suitable for conjugate gradients, we introduce the discrete operator
\begin{equation}\label{eq:poisson_discrete_operator}
(Ay)_{i,j} := (-\Delta_h y)_{i,j}
= \frac{4y_{i,j} - \bigl(y_{i+1,j}+y_{i-1,j}+y_{i,j+1}+y_{i,j-1}\bigr)}{h^2},
\qquad 1\le i,j \le N-2 .
\end{equation}
Homogeneous Dirichlet boundary conditions are imposed strongly by setting $y_{i,j}=0$ for $(i,j)\in\partial\Omega_h$ (analogously for $p$).
The discrete state and adjoint equations are then solved on the interior grid as
\begin{equation}\label{eq:poisson_discrete_state}
A y = u, 
\end{equation}
and
\begin{equation}\label{eq:poisson_discrete_adjoint}
A p = -\frac{\partial}{\partial y}\Bigl(h^2\sum_{i,j}(y_{i,j}-y_{d,i,j})^2\Bigr)
= -2h^2\,(y-y_d),
\end{equation}
We solve \eqref{eq:poisson_discrete_state} and \eqref{eq:poisson_discrete_adjoint} with the conjugate gradient method (CG) provided by JAX, using the same discrete operator $A$ in both solves. The CG was run to a tolerance of $10^{-8}$ with a maximum of 2000 iterations. 
\subsection{Nonlinear Transport Control: Viscous Burgers' Equation}

For the adjoint method, the state and adjoint equations for the viscous Burgers tracking problem
\eqref{eq:burgers} are obtained from the first-order optimality conditions. We recall the state equation

\begin{equation}\label{eq_app:burgers_state}
\begin{aligned}
\frac{\partial y(x,t)}{\partial t}
\;+\; y(x,t)\,\frac{\partial y(x,t)}{\partial x}
&= 0.01\,\frac{\partial^2 y(x,t)}{\partial x^2} \;+\; u(x,t),
\quad (x,t)\in\Omega\times(0,T],\\
y(x,0) &= 0, \qquad x\in\Omega,\\
y(x,t) &= 0, \qquad x\in\partial\Omega,\ \ t\in[0,T],
\end{aligned}
\end{equation}
from which the adjoint equation can be derived as
\begin{equation}\label{eq_app:burgers_adjoint}
\begin{aligned}
-\frac{\partial p(x,t)}{\partial t}
\;-\; y(x,t)\,\frac{\partial p(x,t)}{\partial x}
\;-\; 0.01\,\frac{\partial^2 p(x,t)}{\partial x^2}
&= 2\bigl(y(x,t)-y_d(x,t)\bigr),
\quad (x,t)\in\Omega\times[0,T),\\
p(x,T) &= 0, \qquad x\in\Omega,\\
p(x,t) &= 0, \qquad x\in\partial\Omega,\ \ t\in[0,T].
\end{aligned}
\end{equation}

Given a control $u_k$, one first solves \eqref{eq_app:burgers_state} forward in time for $y_k$ and then solves \eqref{eq_app:burgers_adjoint} backward in time for $p_k$.

To update the control, we compute the gradient of the Lagrangian with respect to $u$.
For an $L_2$ regularization term $\lambda\|u\|_{2}^2$ on $\Omega\times(0,T)$, the gradient is
\begin{equation}\label{eq_app:burgers_adjoint_L2_grad}
\nabla_u \mathcal{L}(u,y,p)=2\lambda u(x,t) - p(x,t).
\end{equation}
For an $H^1$-seminorm regularization term $\lambda\|\nabla u\|_{2}^2$, the gradient takes the form
\begin{equation}\label{eq_app:burgers_adjoint_H1_grad}
\nabla_u \mathcal{L}(u,y,p)= -2\lambda\bigl(\frac{\partial u(x,t)^2}{\partial x^2}+\frac{\partial u(x,t)^2}{\partial t^2}\bigr) - p(x,t),
\end{equation}
with corresponding natural boundary conditions for $u$ on $\partial\Omega$ or at $t\in\{0,T\}$.

We perform a gradient descent step,
\begin{equation}\label{eq_app:burgers_gradient_descent}
u_{k+1}(x,t) = u_k(x,t) - \gamma\,\nabla_u \mathcal{L}(u_k,y_k,p_k),
\end{equation}
where $k$ denotes the iteration index and $\gamma>0$ is the step size. In our implementation, each
iteration consists of:
\begin{enumerate}
    \item Solve the state $y_k$ from \eqref{eq_app:burgers_state} forward in time for the current control $u_k$.
    \item Solve the adjoint $p_k$ from \eqref{eq_app:burgers_adjoint} backward in time using the state $y_k$.
    \item Evaluate the gradient via \eqref{eq_app:burgers_adjoint_L2_grad} or \eqref{eq_app:burgers_adjoint_H1_grad}.
    \item Update the control using \eqref{eq_app:burgers_gradient_descent}.
\end{enumerate}

To solve the state equation \eqref{eq_app:burgers_state} and its adjoint \eqref{eq_app:burgers_adjoint},
we discretize $\Omega\times[0,T]$ with a uniform grid
$x_i = ih$, $i=0,\dots,N_x-1$, and $t^n = n\Delta t$, $n=0,\dots,N_t-1$, where
$h = 1/(N_x-1)$ and $\Delta t = T/(N_t-1)$. We use the same sensor grid as the neural operator, i.e., a $(32,32)$ grid. Homogeneous Dirichlet boundary conditions are imposed
strongly by setting $y_0^n=y_{N_x-1}^n=0$ (and analogously for the adjoint).

We advance the viscous Burgers dynamics with an IMEX scheme \cite{AscherRuuthWetton1995}: convection is treated explicitly via a
Rusanov (local Lax--Friedrichs) flux in conservative form, while diffusion is treated implicitly
with a second-order finite difference. Writing $f(y)=\tfrac12 y^2$, the numerical flux at $x_{i+\frac12}$
is
\begin{equation}\label{eq_app:burgers_rusanov_flux}
F_{i+\frac12}^n
=\frac12\Bigl(f(y_i^n)+f(y_{i+1}^n)\Bigr)
-\frac{\alpha}{2}\max\bigl(|y_i^n|,|y_{i+1}^n|\bigr)\bigl(y_{i+1}^n-y_i^n\bigr),
\end{equation}
where $\alpha>0$ is a stabilization parameter. The discrete convection term on interior nodes
$i=1,\dots,N_x-2$ is then
\begin{equation}\label{eq_app:burgers_conv_disc}
\bigl(\mathcal{C}(y^n)\bigr)_i \;=\; -\frac{F_{i+\frac12}^n - F_{i-\frac12}^n}{h}.
\end{equation}
For diffusion, we use the centered Laplacian
\begin{equation}\label{eq_app:burgers_lap_disc}
(D_{xx} y^{n+1})_i \;=\; \frac{y_{i+1}^{n+1}-2y_i^{n+1}+y_{i-1}^{n+1}}{h^2}.
\end{equation}
The full-time step reads
\begin{equation}\label{eq_app:burgers_imex_step}
\Bigl(I - \Delta t\,\nu D_{xx}\Bigr)y^{n+1}
=
y^n + \Delta t\Bigl(\mathcal{C}(y^n) + u^n\Bigr),
\end{equation}
which yields, on interior indices, a tridiagonal system solved by a Thomas algorithm.

The tracking objective is discretized with the grid-weighted $L_2$ inner product,
\begin{equation}\label{eq_app:burgers_misfit_disc}
\int_0^T\!\!\int_\Omega (y-y_d)^2\,dx\,dt
\;\approx\;
\Delta t\,h \sum_{n,i}\bigl(y_i^n-y_{d,i}^n\bigr)^2,
\end{equation}
and we use the same weighting for the control regularization,
e.g. $\|u\|_{L_2}^2 \approx \Delta t\,h\sum_{n,i}(u_i^n)^2$ and
$\|\partial_x u\|_{L_2}^2 \approx \Delta t\,h\sum_{n,i}\bigl((u_{i+1}^n-u_i^n)/h\bigr)^2$.
After discretization, the full (discrete) adjoint gradient is evaluated by a backward-in-time sweep
using the transpose Jacobian of the one-step update \eqref{eq_app:burgers_imex_step}.

\subsection{Stokes Flow}

We generate reference solutions for the steady incompressible Stokes equations on
$\Omega=[0,1]\times[0,1]$ with a forcing (control) field $\mathbf{u}=(u_1,u_2)$:
\begin{equation}\label{eq_app:stokes}
\begin{aligned}
-\nu \Delta \mathbf{v}(\mathbf{x}) + \nabla p(\mathbf{x}) &= \mathbf{u}(\mathbf{x}),
\qquad \mathbf{x}\in\Omega,\\
\nabla\cdot \mathbf{v}(\mathbf{x}) &= 0, \qquad \mathbf{x}\in\Omega,
\end{aligned}
\end{equation}
where $\mathbf{v}=(v_1,v_2)$ is the velocity and $p$ is the pressure. We impose a parabolic inflow at
$\Gamma_{\mathrm{in}}=\{0\}\times[0,1]$,
\begin{equation}\label{eq_app:stokes_inlet}
v_1(0,y)=a\,4U_{\max}y(1-y),\qquad v_2(0,y)=0,
\end{equation}
no-slip walls $\Gamma_{\mathrm{w}}=[0,1]\times\{0,1\}$ with $\mathbf{v}=\mathbf{0}$, and an outflow
condition $\partial_x v_1=\partial_x v_2=0$ at $\Gamma_{\mathrm{out}}=\{1\}\times[0,1]$. To remove the
pressure nullspace we enforce a gauge condition by subtracting the mean, $\int_\Omega p\,d\mathbf{x}=0$.

We discretize $\Omega$ on a uniform node grid $\{(x_i,y_j)\}$ with
$x_i = i\Delta x$, $i=0,\dots,N_x-1$, $y_j=j\Delta y$, $j=0,\dots,N_y-1$, where
$\Delta x=1/(N_x-1)$ and $\Delta y=1/(N_y-1)$. Spatial derivatives are approximated by second-order
finite differences on nodes: for a scalar field $\phi$,
\begin{equation}\label{eq_app:stokes_fd_ops}
\frac{\partial \phi}{\partial x}\Big|_{i,j} \approx \frac{\phi_{i+1,j}-\phi_{i-1,j}}{2\Delta x},
\qquad
\frac{\partial \phi}{\partial y}\Big|_{i,j} \approx \frac{\phi_{i,j+1}-\phi_{i,j-1}}{2\Delta y},
\end{equation}
and
\begin{equation}\label{eq_app:stokes_fd_lap}
\Delta \phi\Big|_{i,j} \approx
\frac{\phi_{i+1,j}-2\phi_{i,j}+\phi_{i-1,j}}{\Delta x^2}
+\frac{\phi_{i,j+1}-2\phi_{i,j}+\phi_{i,j-1}}{\Delta y^2}.
\end{equation}
Dirichlet conditions on the inlet and walls are imposed strongly by overwriting boundary nodes, while the outflow Neumann condition is enforced via a ghost-node relation (equivalently, copying from the last interior column).

We solve \eqref{eq_app:stokes} with a pressure-correction (projection) fixed-point iteration. Given a
current pressure $p^{(k)}$, we compute an intermediate velocity $\mathbf{v}^\star$ from two decoupled
Poisson problems,
\begin{equation}\label{eq_app:stokes_momentum_disc}
-\nu \Delta v_1^\star + \frac{\partial p^{(k)}}{\partial x} = u_1,
\qquad
-\nu \Delta v_2^\star + \frac{\partial p^{(k)}}{\partial y} = u_2,
\end{equation}
subject to the velocity boundary conditions. Next, we solve a Poisson equation for the pressure
correction $\phi$,
\begin{equation}\label{eq_app:stokes_phi}
\Delta \phi = \frac{1}{\alpha}\,\nabla\cdot \mathbf{v}^\star,
\end{equation}
with homogeneous Neumann boundary conditions on all sides and a single pinned value (e.g.,
$\phi(0,0)=0$) to remove the constant nullspace. We then update
\begin{equation}\label{eq_app:stokes_correction}
\mathbf{v}^{(k+1)} = \mathbf{v}^\star - \alpha \nabla \phi,
\qquad
p^{(k+1)} = p^{(k)} + \phi,
\end{equation}
and iterate until the maximum interior divergence and the maximum velocity update fall below
prescribed tolerances. Each Poisson subproblem is solved with the conjugate gradient method (CG).

\subsection{Regularization Parameter Sweep}\label{appendix_b_regularization}
Table \ref{tab:ref_summary_control_mse_transposed} shows the control error (MSE) for different regularizations for the scalar elliptic and the nonlinear transport control problems. From this table best regularization values were selected for comparison in section \ref{section_5}.
\begin{table}[h]
\centering
\scriptsize
\caption{Control error (MSE) for $L_2$ and $H^1$ regularization across the scalar elliptic control and the nonlinear transport control problems for different regularization parameters. For each problem, the best result is bolded separately within the $L_2$ block and within the $H^1$ block.}
\label{tab:ref_summary_control_mse_transposed}
\renewcommand{\arraystretch}{1.15}
\setlength{\tabcolsep}{4pt}

\begin{tabular}{c c cccc}
\toprule
Method & Regularization & P1 & P2 & $J_{Burg,1}$ & $J_{Burg,2}$ \\
\midrule
\multirow{9}{*}{$L_2$} & $\lambda=10^{-2}$ & $1.48\times 10^{-1}$ & $2.47\times 10^{-1}$ & $2.17\times 10^{-1}$ & $1.99\times 10^{-1}$ \\
 & $\lambda=10^{-3}$ & $1.85\times 10^{-2}$ & $1.76\times 10^{-1}$ & $\bm{1.72\times 10^{-1}}$ & $1.07\times 10^{-1}$ \\
 & $\lambda=10^{-4}$ & $3.43\times 10^{-4}$ & $7.06\times 10^{-2}$ & $2.58\times 10^{-1}$ & $\bm{9.84\times 10^{-2}}$ \\
 & $\lambda=10^{-5}$ & $\bm{5.26\times 10^{-6}}$ & $3.47\times 10^{-2}$ & $2.84\times 10^{-1}$ & $9.87\times 10^{-2}$ \\
 & $\lambda=10^{-6}$ & $1.45\times 10^{-5}$ & $2.48\times 10^{-2}$ & $2.88\times 10^{-1}$ & $9.88\times 10^{-2}$ \\
 & $\lambda=10^{-7}$ & $1.1\times 10^{-4}$ & $2.2\times 10^{-2}$ & $2.88\times 10^{-1}$ & $9.88\times 10^{-2}$ \\
 & $\lambda=10^{-8}$ & $1.35\times 10^{-4}$ & $\bm{2.18\times 10^{-2}}$ & $2.88\times 10^{-1}$ & $9.88\times 10^{-2}$ \\
 & $\lambda=10^{-9}$ & $1.38\times 10^{-4}$ & $2.18\times 10^{-2}$ & $2.88\times 10^{-1}$ & $9.88\times 10^{-2}$ \\
 & $\lambda=10^{-10}$ & $1.38\times 10^{-4}$ & $2.18\times 10^{-2}$ & $2.88\times 10^{-1}$ & $9.88\times 10^{-2}$ \\
\midrule
\multirow{9}{*}{$H^1$} & $\lambda=10^{-2}$ & $2.28\times 10^{-1}$ & $2.6\times 10^{-1}$ & $6.86\times 10^{-1}$ & $4.48\times 10^{-1}$ \\
 & $\lambda=10^{-3}$ & $1.78\times 10^{-1}$ & $2.57\times 10^{-1}$ & $2.83\times 10^{-1}$ & $3.49\times 10^{-1}$ \\
 & $\lambda=10^{-4}$ & $3.92\times 10^{-2}$ & $2.34\times 10^{-1}$ & $1.24\times 10^{-1}$ & $2.07\times 10^{-1}$ \\
 & $\lambda=10^{-5}$ & $2.15\times 10^{-3}$ & $1.39\times 10^{-1}$ & $\bm{1.04\times 10^{-1}}$ & $9.32\times 10^{-2}$ \\
 & $\lambda=10^{-6}$ & $9.25\times 10^{-4}$ & $6.27\times 10^{-2}$ & $1.55\times 10^{-1}$ & $6.28\times 10^{-2}$ \\
 & $\lambda=10^{-7}$ & $6.08\times 10^{-4}$ & $2.93\times 10^{-2}$ & $2.48\times 10^{-1}$ & $5.73\times 10^{-2}$ \\
 & $\lambda=10^{-8}$ & $3.95\times 10^{-4}$ & $1.68\times 10^{-2}$ & $2.83\times 10^{-1}$ & $5.66\times 10^{-2}$ \\
 & $\lambda=10^{-9}$ & $\bm{3.62\times 10^{-4}}$ & $\bm{1.05\times 10^{-2}}$ & $2.84\times 10^{-1}$ & $\bm{5.65\times 10^{-2}}$ \\
 & $\lambda=10^{-10}$ & $4.68\times 10^{-4}$ & $1.13\times 10^{-2}$ & $2.67\times 10^{-1}$ & $6.73\times 10^{-2}$ \\
\bottomrule
\end{tabular}
\end{table}

\newpage
\section{Additional ODE-constrained Optimal Control Problems}\label{appendix_d}
In this appendix, we do a preliminary study of the method for additional ODE problems, not of the tracking type, but of quadratic and terminal cost. In contrast to the tracking problems presented in the main text, these models require an additional smoothing regularization to converge to the solution. We include this brief study to demonstrate that our method of using a physics-informed neural operator as a surrogate model is not limited to tracking problems or reachable controls.

\subsection{Nonlinear ODE}\label{section_nonlinear_ode}
This dynamical system has been used as a benchmark by \cite{Rentsen2023Computational} and \cite{Garg2010Unified} for solving optimal control problems. The dynamical system is described by the following equations.
  \begin{equation}\label{eq:nonlinear_ode}
  \begin{aligned}
    \frac{dy(t)}{dt} &= \tfrac{5}{2}\Big(-y(t) + y(t)u(t) - u(t)^2\Big),\\
    y(0) &= 1,\quad t \in [0,1].
  \end{aligned}
  \end{equation}
We choose to start with this easy problem to demonstrate that the method works for both terminal cost functions and quadratic cost functions, hence
\begin{subequations}\label{cost:ode}
  \begin{align}
    J_{\mathrm{ODE},1}(u)
      &= -y(1)\\
    J_{\mathrm{ODE},2}(u)
      &= \int_{0}^{1}(y(t)^2+u(t)^2)dt \label{cost:ode:b}
  \end{align}
  \end{subequations}
We train a physics-informed DeepONet with a similar setup as demonstrated in \ref{appendix_a}. For the branch and trunk network, we use the modified network described in \cite{wang2021learning}, with symmetric sizes of 5 hidden layers and 300 neurons. We use the hyperbolic tangent activation function. The sensor grid was discretized to 100 points. 

The training data consisted of 300000 functions, where we constructed the dataset similarly to that for the PDE-tracking problems. We used GRF functions with a lengthscale in $[0.1,1.0]$ as well as random polynomials of up to degree 5. We projected the input functions to remain in a domain of $[-1.2,1.2]$. For other hyperparameters of training, we kept them similar as for the PDE-tracking models. 

For the optimizing we used a penalty parameter $\mu=50$ and a regularization parameter $\lambda=1$. In contrast to PDE-tracking, we used a $H^2$-seminorm, i.e., $\|\Delta u\|_2^2$ regularization in addition to the residual penalty for this problem. We noticed that the $H^2$-seminorm provided smoother results, but the problem is also solvable with the $H^1$-seminorm regularization. We used the AdamW optimizer of Optax \cite{deepmind2020jax} with update step $\gamma=0.001$ without any step scheduler. We ran the optimization for 2000 iteration steps.

For the reference methods, we used the CasADi \cite{Andersson2018} library with the same grid size as the neural operator. The problem was solved with a single shooting method, where an explicit 4th order Runge-Kutta integrator was used for time-integration. Finally, we used IPOPT for the background NLP solver within CasADi.

\begin{figure}[t]  
    \centering
    \includegraphics[width=0.99\linewidth]{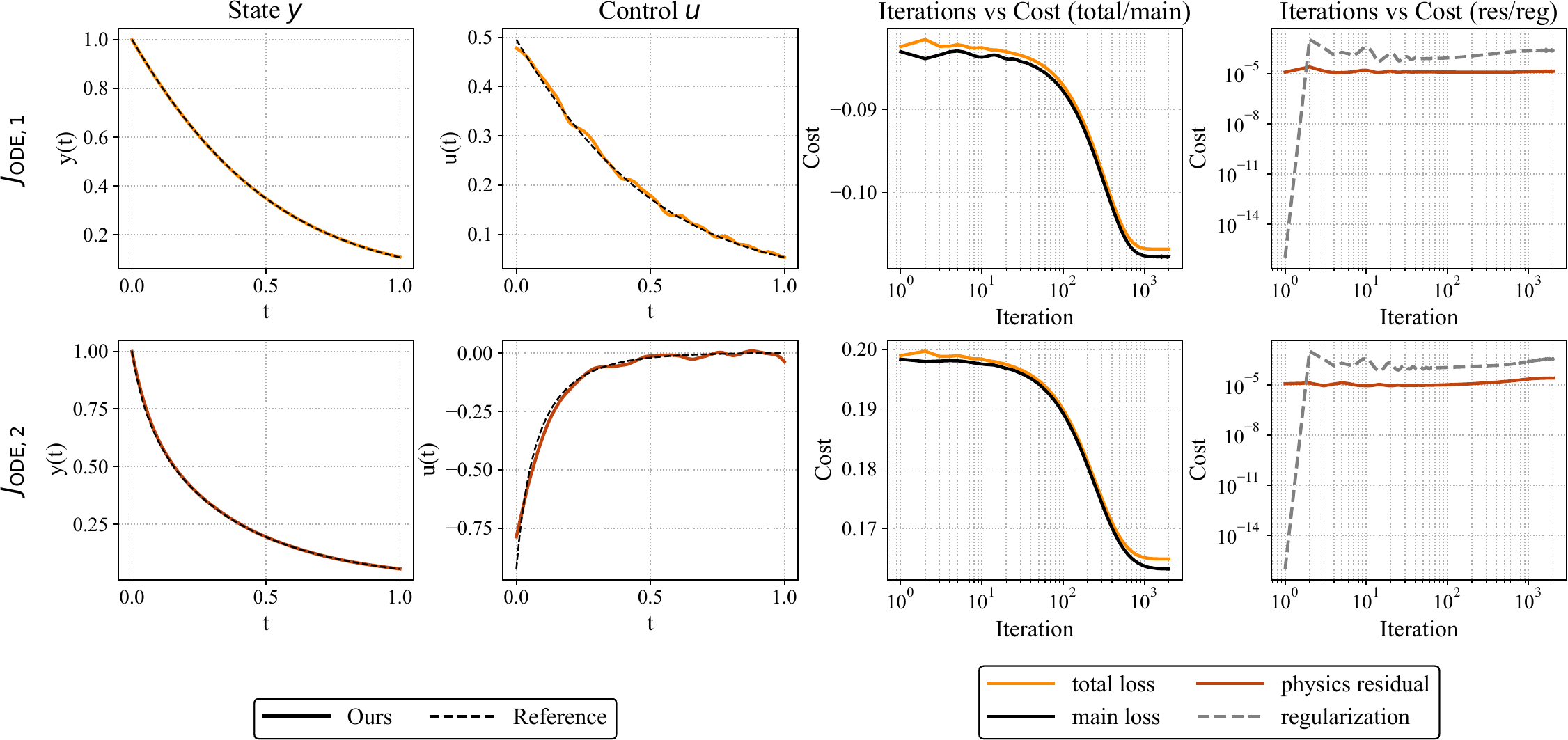}  
    \caption{Nonlinear ODE control problems \eqref{eq:nonlinear_ode}: comparison of the optimized solution (solid lines) against the reference solution (dashed lines) for Problem 1 (top row) and Problem 2 (bottom row). Columns (left to right) show: state trajectory $y(t)$, control function $u(t)$, cost convergence (total vs. main loss), and residual convergence (physics residual vs. regularization).}
    \label{fig:nonlinear_ode_appendix}
\end{figure}

The results of the optimization and the reference are shown in Figure \ref{fig:nonlinear_ode_appendix}. The optimizer finds an approximate solution that aligns with the reference solution. The residual remains approximately constant during the optimization process. The jump in the regularization term is due to the cold start of the control, which is initialized as $u(t)=0$, and thus the $H^2$-seminorm of a zero vector is zero.

As a conclusion, a physics-informed neural operator can act as a surrogate model for the dynamics. We showed that our method worked for terminal cost problems and quadratic minimization problems (quadratic cost), but it should be regularized with a smoothing regularizer, such as $H^1$ or $H^2$. However, as we do not have a theoretical explanation for this, we continue to study the limitations and requirements for our method

\subsection{SIR-model with contact and vaccination control}
This epidemic control model follows the standard Susceptible-Infected-Recovered (SIR) formulation with an additional vaccination control $v(t)$, similar to the setting in~\cite{Marinov2022}, and a control of contacts, $u(t)$. The compartmental models are widely studied as an application of optimal control in epidemiology. We choose this problem to demonstrate that the method works for multidimensional ODEs and for multiple control functions. The dynamics are given by

  \begin{equation}\label{eq:sir}
  \begin{aligned}
    \frac{dS(t)}{dt} &= -\beta(1-u(t))S(t)I(t) - v(t)S(t),\\
    \frac{dI(t)}{dt} &= \beta(1-u(t))S(t)I(t) - \gamma I(t),\\
    \frac{dR(t)}{dt} &= \gamma I(t) + v(t)S(t),\\
    \text{where }
    & S(0) = 0.98, I(0) = 0.02,,R(0) = 0,\\ 
    \beta=4,\, &\gamma=0.5,\; t \in [0,10].
  \end{aligned}
  \end{equation}

We choose two cost functions, the first one considers vaccination and control planning for minimizing costs. The second incorporates healthcare capacity constraints by adding a penalty on the infected population whenever it exceeds 20\%,
  \begin{subequations}\label{cost:sir}
  \begin{align}
    J_{\mathrm{SIR},1}(u,v)
      &= \int_0^{10} \!\Big( I^2 + 0.2 u^2 + 5 v^2 \Big)\,dt, \label{cost:sir:a}\\
    J_{\mathrm{SIR},2}(u,v)
    &= \int_0^{10} \!\Big(0.1 u^2 + v^2 + 10\bigl[I-0.2\bigr]_+^2 \Big)\,dt, \label{cost:sir:b}
  \end{align}
  \end{subequations}

We train a physics-informed DeepONet with a similar setup as demonstrated in \ref{appendix_a}. For the branch and trunk network, we use the modified network described in \cite{wang2021learning}, with symmetric sizes of 4 hidden layers and 600 neurons. We use the hyperbolic tangent activation function. The sensor grid was discretized to 200 points. We used to separate branch networks for the controls $u(t)$ and $v(t)$ and combined them with a dot product after the last layer.

The training data consisted of 300000 functions, where we constructed the dataset similarly to that for the PDE-tracking problems. We used GRF functions with a lengthscale in $[0.05,0.5]$ as well as random polynomials of up to degree 5. We projected the input functions to remain in a domain of $[0,1.0]$. For other hyperparameters of training, we kept them similar as for the PDE-tracking models. 

For the optimizing we used a penalty parameter $\mu=100$ and a regularization parameter $\lambda=0.01$. In contrast to PDE-tracking, we used a $H^2$-seminorm, i.e., $\|\Delta u\|_2^2$ regularization in addition to the residual penalty for this problem. 

For this problem, we also used CasADi \cite{Andersson2018} with the same settings as for the nonlinear ODE in section \ref{section_nonlinear_ode}. We again used the same grid size as the DeepONet model, i.e., a grid size of 200.

\begin{figure}[t]  
    \centering
    \includegraphics[width=0.99\linewidth]{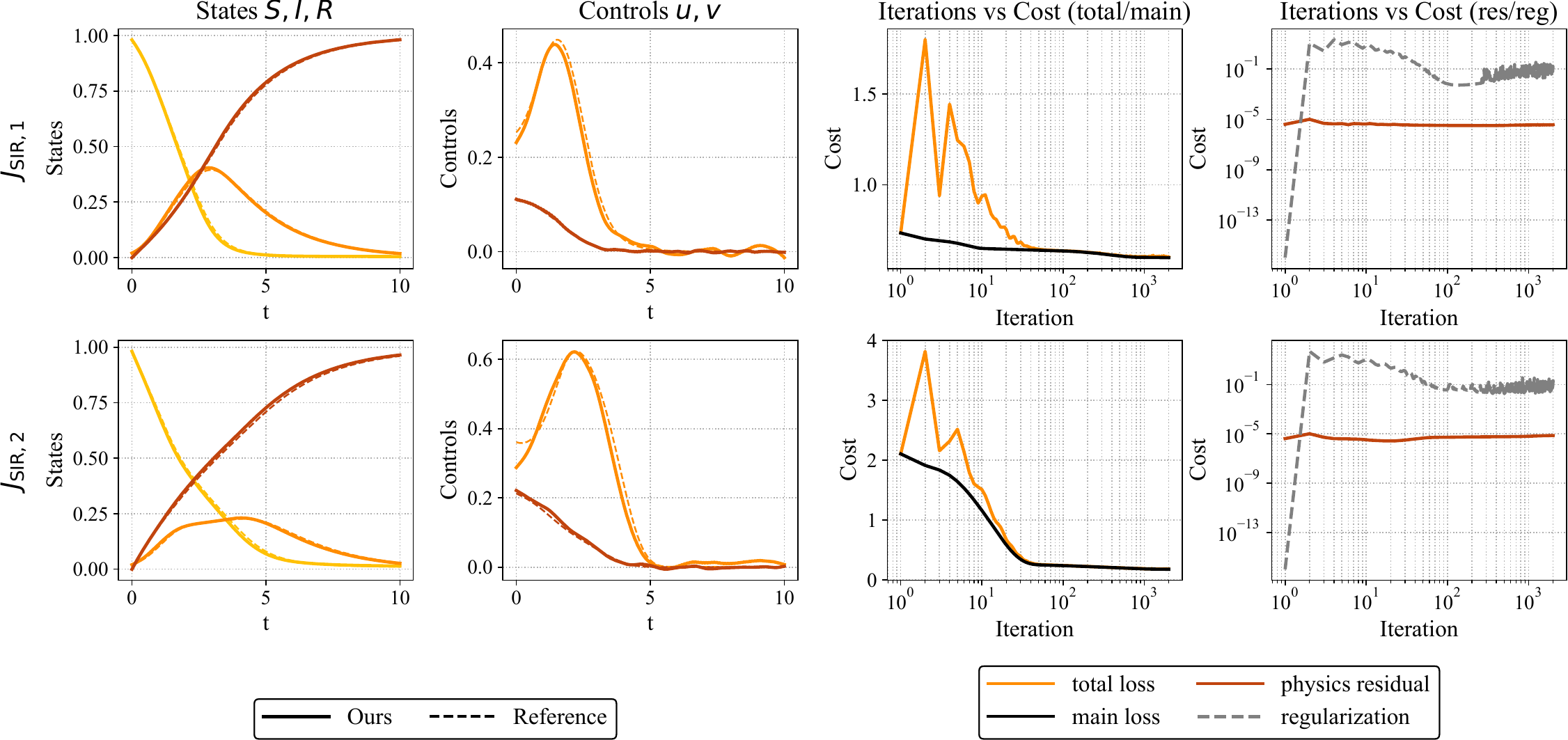}  
    \caption{SIR control problems \eqref{eq:sir}: comparison of the optimized solution (solid lines) against the reference solution (dashed lines) for Problem 1 (top row) and Problem 2 (bottom row). Columns (left to right) show: state trajectories $S, I, R$, control functions $u, v$, cost convergence (total vs. main loss), and residual convergence (physics residual vs. regularization).}
    \label{fig:sir_appendix}
\end{figure}

As shown in Figure \ref{fig:sir_appendix}, our method is capable of finding the solution to the control problem \eqref{cost:sir}. Similarly to the nonlinear ODE, the residual remains approximately constant during the iterations. Thus, the physics-informed neural operator is capable of acting as a surrogate model for a multidimensional ODE with multiple controls.

\newpage
\bibliographystyle{elsarticle-num} 

\bibliography{references}

\end{document}